\documentclass[12pt,twoside,a4paper]{article}
\usepackage{mathrsfs}
\usepackage{mathrsfs}
\usepackage{mathrsfs}
\usepackage{mathrsfs}
\usepackage{mathrsfs}
\usepackage{mathrsfs}
\usepackage{mathrsfs}
\usepackage{stmaryrd}
\usepackage{amsfonts}
\usepackage{amsmath,amssymb}
\usepackage{mathrsfs}
\usepackage{hyperref}
\usepackage{graphicx}
\usepackage{graphicx}

\newtheorem{thm}{Theorem}[section]
\newtheorem{cor}[thm]{Corollary}
\newtheorem{lem}[thm]{Lemma}
\newtheorem{prop}[thm]{Proposition}

\newtheorem{rem}[thm]{Remark}
\newtheorem{con}[thm]{Conjecture}
\numberwithin{equation}{section}\allowdisplaybreaks

\def\leq{\leqslant}
\def\ge{\geqslant}
\def\leq{\leqslant}
\def\geq{\geqslant}

\begin{document}

\title{\large\bf Global Smooth Effects, Well-Posedness and Scattering \\ for the Derivative Nonlinear Schr\"odinger  Equation \\ with Small Rough Data}

\author{\normalsize \bf Wang Baoxiang, \ Han Lijia, \  Huang Chunyan\date{}\\
{\footnotesize \it
LMAM, School of Mathematical Sciences, Peking University, Beijing 100871, China}\\
{\small E-mails: wbx,\ hanlijia, hcy@math.pku.edu.cn } }\maketitle


\begin{abstract}
\rm  We obtain the global smooth effects  for the solutions of the
linear Schr\"odinger equation in anisotropic Lebesgue spaces.
Applying these estimates, we study the Cauchy problem for the
generalized elliptical and non-elliptical derivative nonlinear
Schr\"odinger equations (DNLS) and get the global well posedness of
solutions with small data  in modulation spaces
$M^{3/2}_{2,1}(\mathbb{R}^n)$. Noticing that $B^{s+n/2}_{2,1}
\subset M^s_{2,1} \subset B^s_{2,1}$ are optimal inclusions, we have
shown the global well posedness of DNLS with a class of rough data.  As by products,
the existence of the scattering operators with small data is also obtained.\\

\\

{\it MSC:} 35 Q 55, 46 E 35, 47 D 08.
\end{abstract}

\section{Introduction}
This paper is a continuation of our earlier work \cite{WaWa} and we
study the Cauchy problem for the generalized derivative nonlinear
Schr\"odinger equation (gDNLS)
\begin{align}
{\rm i} u_t + \Delta_\pm u = F(u, \bar{u}, \nabla u, \nabla
\bar{u}), \quad u(0,x)= u_0(x),
 \label{gNLS}
\end{align}
where $u$ is a complex valued function of $(t,x)\in \mathbb{R}
\times \mathbb{R}^n$,
\begin{align}
\Delta_\pm u = \sum^n_{i=1} \varepsilon_i \partial^2_{x_i}, \quad
\varepsilon_i \in \{1,\,  -1\}, \quad i=1,...,n,
 \label{Delta-pm}
\end{align}
$\nabla =(\partial_{x_1},..., \partial_{x_n})$, $F:
\mathbb{C}^{2n+2} \to \mathbb{C}$ is a polynomial series,
\begin{align}
F(z) = F(z_1,..., z_{2n+2})= \sum^{}_{m+1<|\beta|<\infty} c_\beta
z^\beta, \quad c_\beta \in \mathbb{C},
 \label{poly}
\end{align}
$2 \le m <\infty$, $m\in \mathbb{N}$, $\sup_\beta
|c_\beta|<\infty$\footnote{In fact, $c_\beta$ is not necessarily
bounded, condition $\sup_\beta |c_\beta|<\infty$ can be replaced by
$|c_\beta |\le C^{|\beta|}$.}. A typical nonlinear term is the
following
$$
F(u, \bar{u}, \nabla u, \nabla \bar{u}) = |u|^2\vec{\lambda}\cdot
\nabla u + u^2\vec{\mu}\cdot \nabla \bar{u} +|u|^2 u,
$$
which is a model equation in the strongly interacting many-body
systems near criticality as recently described in terms of nonlinear
dynamics \cite{TD, DT, CT}. Another typical nonlinearity is
$$
F(u, \bar{u}, \nabla u, \nabla \bar{u})= (1\mp |u|^2)^{-1}|\nabla
u|^2 u= \sum^\infty_{k=0} \pm |u|^{2k} |\nabla u|^2 u,  \quad |u|<1,
$$
which is a deformation of the Schr\"odinger map equation \cite{DW,
Io-Ke}.

A large amount of work has been devoted to the study of the local
and global well posedness of \eqref{gNLS}, see Bejenaru and Tataru
\cite{Be-Ta}, Chihara \cite{Chih1,Chih2}, Kenig, Ponce and Vega
\cite{KePoVe1,KePoVe2}, Klainerman \cite{Klai}, Klainerman and Ponce
\cite{Kl-Po}, Ozawa and Zhai \cite{Oz-Zh}, Shatah \cite{Shat}, B.
Wang and Y. Wang \cite{WaWa}. When the nonlinear term $F$ satisfies
some energy structure conditions, or the initial data suitably
decay, the energy method, which went back to the work of Klainerman
\cite{Klai} and was developed in \cite{Chih1,Chih2,Kl-Po,Oz-Zh,
Shat}, yields the global existence of \eqref{gNLS} in the elliptical
case $\Delta_\pm= \Delta$. Recently, Ozawa and Zhai obtained the
global well posedness in $H^{s}(\mathbb{R}^n)$ ($n\ge 3$, $s>2+n/2$,
$m\ge 2$) with small data for \eqref{gNLS} in the elliptical case,
where an energy structure condition on $F$ is still required.

By setting up the local smooth effects for the solutions of the
linear Schr\"odinger equation, Kenig, Ponce and Vega
\cite{KePoVe1,KePoVe2} were able to deal with the non-elliptical
case and they established the local well posedness of Eq.
\eqref{gNLS} in $H^s$ with $s\gg n/2$.  Recently, the local well
posedness results have been generalized to the quasi-linear
(ultrahyperbolic) Schr\"odinger equations, see
\cite{KePoVe3,KePoRoVe}.

In one spatial dimension, B. Wang and Y. Wang \cite{WaWa} showed the
global well posedness of gDNLS \eqref{gNLS} for small data in
critical Besov spaces $\dot B^{1+ n/2-2/m}_{2,1}\cap $ $\dot B^{1+
n/2-1/M}_{2,1}(\mathbb{R})$, $m\ge 4$. In higher spatial dimensions
$n\ge 2$, by using Kenig, Ponce and Vega's local smooth effects and
establishing time-global maximal function estimates in space-local
Lebesgue spaces, B. Wang and Y. Wang \cite{WaWa} showed the global
well posedness of gDNLS \eqref{gNLS} for small data in Besov spaces
$B^s_{2,1}(\mathbb{R}^n)$ with $s>n/2+3/2$, $m\ge 2+4/n$.

Wang and Huang \cite{WaHu} obtained the global well posedness of
\eqref{gNLS} in one spatial dimension with initial data in
$M^{1+1/m}_{2,1}$, $m\ge 4$. In this paper, we will use a new way to
study the global well posedness of \eqref{gNLS} and show that
\eqref{gNLS} is globally well posed in $M^{s}_{2,1}(\mathbb{R}^n)$
with $s\ge 3/2, \,\, m\ge 2$ and $m>4/n$ for the small Cauchy data.
Our starting point is the smooth effect estimates for the linear
Schr\"odinger equation in one spatial dimension (cf.
\cite{CS,KePoVe,KePoVe1,Sj,Ve}),  from which we get a series of
linear estimates in higher dimensional anisotropic Lebesgue spaces,
including the global smooth effect estimates, the maximal function
estimates and their relations to the Strichartz estimates. The
maximal function estimates follows an idea as in Ionescu and Kenig
\cite{Io-Ke}.  These estimates together with the frequency-uniform
decomposition method yield the global well posedness and scattering
of solutions in modulation spaces $M^s_{2,1}$, $s\ge 3/2$.

\subsection{$M^s_{2,1}$ and $B^s_{2,1}$}

In this paper, we are mainly interested in the cases that the
initial data $u_0$ belongs to the modulation space $ M^s_{2,1}$ for
which the norm can be equivalently defined in the following way (cf.
\cite{Fei2,Wa1,WaHe,WaHu}):
\begin{align}
\|f\|_{M^s_{2,1}}= \sum_{k\in \mathbb{Z}^n} \langle k \rangle^s
\|\mathscr{F} f\|_{L^2(Q_k)},  \label{mod-s21}
\end{align}
where $\langle k\rangle=1+|k|$, $Q_k= \{\xi: -1/2\le \xi_i-k_i <1/2,
\ i=1,...,n\}$. For simplicity, we write $M_{2,1}=M^0_{2,1}$.  Since
only the modulation space $M^s_{2,1}$ will be used in this paper, we
will not state the defination of the general modulation spaces
$M^s_{p,q}$, one can refer to Feichtinger \cite{Fei2}. Modulation
spaces $M^s_{2,1}$ are related to the Besov spaces $B^s_{2,1}$ for
which the norm is defined as follows:
\begin{align}
\|f\|_{B^s_{2,1}}=\|\mathscr{F} f\|_{L^2(B(0,1))}+ \sum^\infty_{j=
1} 2^{sj} \|\mathscr{F} f\|_{L^2(B(0,2^j) \setminus B(0, 2^{j-1}))},
\label{Besov-s21}
\end{align}
where $B(x_0,R):= \{\xi\in \mathbb{R}: \ |\xi-x_0|\le R \}$.  It is
known that there holds the following  optimal inclusions between $
B^{n/2+s}_{2,1}$, $ M^s_{2,1}$ and $B^s_{2,1}$ (cf.
\cite{Toft,SuTo,WaHu}):
\begin{align}
B^{n/2+s}_{2,1}\subset M^s_{2,1} \subset B^s_{2,1}.
\label{Besov-s21}
\end{align}
 So, comparing $M^{s}_{2,1}$ with
$B^{s+n/2}_{2,1}$,  we see that $M^{s}_{2,1}$ contains a class of
functions $u$ satisfying $\|u\|_{M^s_{2,1}}= \infty$ but
$\|u\|_{B^{s+n/2}_{2,1}} \ll 1$. On the other hand, we can also find
a class of rough functions $u$ satisfying $\|u\|_{B^s_{2,1}}=
\infty$ but $\|u\|_{M^{s}_{2,1}} \ll 1$. Another important inclusion
between $M_{2,1}$ and $L^\infty$ is that $M_{2,1}\subset L^\infty$
and this embedding is also optimal, see Figure 1.

\begin{figure}
\begin{center}
\includegraphics[height=6cm,width=6cm]{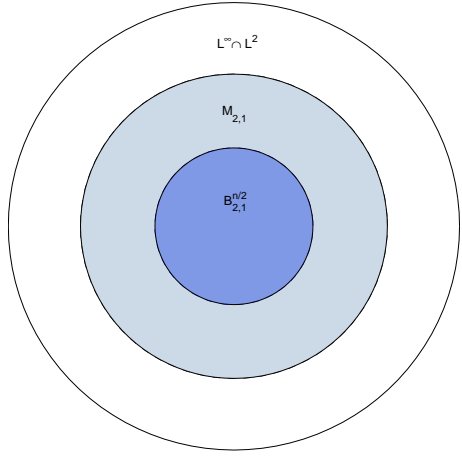}
\begin{minipage}{11cm}
\caption{\small Optimal inclusions: $B^{n/2}_{2,1}\subset M_{2,1}
\subset L^\infty \cap L^2$.}
\end{minipage}
\end{center}
\end{figure}

\subsection{Main Results}

For the definitions of the anisotropic Lebesgue spaces
$L^{p_1}_{x_i} L^{p_2}_{(x_j)_{j\not=
i}}L^{p_2}_t(\mathbb{R}^{1+n})$ and the frequency-uniform
decomposition operators $\{\Box_k\}_{k\in \mathbb{Z}^n}$, one can
refer to Section \ref{Notation}. We have

\begin{thm} \label{DNLS-mod}
Let $n\ge 2$, $ 2< m <\infty, \; m>4/n$. Assume that $u_0 \in
M^{3/2}_{2,1}$ and $\|u_0\|_{M^{3/2}_{2,1}} \le \delta$ for some
small $\delta>0$. Then \eqref{gNLS} has a unique global solution
$u\in C(\mathbb{R}, M^{3/2}_{2,1}) \cap X$, where
\begin{align}
\!\!\!\!\!\! \|u\|_{X} =  &\sum_{\alpha=0,1} \ \sum_{i, \, \ell=1}^n
\ \sum_{k\in \mathbb{Z}^n, \ |k_i|>4} \langle k_i\rangle
\left\|\partial^\alpha_{x_\ell} \Box_k u \right\|_{L^\infty_{x_i}
L^2_{(x_j)_{j\not= i}}L^2_t(\mathbb{R}^{1+n})} \nonumber\\
&  + \sum_{\alpha=0,1} \ \sum_{i,\, \ell=1}^n \ \sum_{k\in
\mathbb{Z}^n} \langle k \rangle^{1/2-1/m} \left\|
\partial^\alpha_{x_\ell} \Box_k u \right\|_{L^{m}_{x_i}
L^\infty_{(x_j)_{j\not= i} }
L^\infty_t(\mathbb{R}^{1+n})} \nonumber\\
& + \sum_{\alpha=0,1} \ \sum_{\ell=1}^n \sum_{k\in \mathbb{Z}^n}
\langle k \rangle^{1/2} \left\|
\partial^\alpha_{x_\ell} \Box_k u \right\|_{L^\infty_t L^2_x  \bigcap  L^{2+m}_{x,t}
(\mathbb{R}^{1+n})}, \label{modX.2}
\end{align}
where $k=(k_1,...,k_n)$ and we further have  $\|u\|_{X} \lesssim
\delta$. Moreover, the scattering operator of \eqref{gNLS} $S$
carries a zero neighborhood in $C(\mathbb{R}, M^{3/2}_{2,1})$ into
$C(\mathbb{R}, M^{3/2}_{2,1})$.
\end{thm}
In Theorem \ref{DNLS-mod}, if $u_0 \in M^s_{2,1}$ with $s>3/2$, then
we have $u\in C(\mathbb{R}, M^{s}_{2,1})$.  If $m=2$, we need to
assume the initial data have stronger regularity:

\begin{thm} \label{DNLS-modm2}
Let $n\ge 3$, $m=2$. Assume that $u_0 \in M^{5/2}_{2,1}$ and
$\|u_0\|_{M^{5/2}_{2,1}} \le \delta$ for some small $\delta>0$. Then
\eqref{gNLS} has a unique global solution $u\in C(\mathbb{R},
M^{5/2}_{2,1}) \cap Y$, where
\begin{align}
\!\!\!\!\!\! \|u\|_{Y} =  &\sum_{\alpha=0,1} \ \sum_{i, \, \ell=1}^n
\ \sum_{k\in \mathbb{Z}^n, \ |k_i|>4} \langle k_i\rangle^{2}
\left\|\partial^\alpha_{x_\ell} \Box_k u \right\|_{L^\infty_{x_i}
L^2_{(x_j)_{j\not= i}}L^2_t(\mathbb{R}^{1+n})} \nonumber\\
&  + \sum_{\alpha=0,1} \ \sum_{i,\, \ell=1}^n \ \sum_{k\in
\mathbb{Z}^n}  \left\|
\partial^\alpha_{x_\ell} \Box_k u \right\|_{L^2_{x_i}
L^\infty_{(x_j)_{j\not= i} }
L^\infty_t(\mathbb{R}^{1+n})} \nonumber\\
& + \sum_{\alpha=0,1} \ \sum_{\ell=1}^n \sum_{k\in \mathbb{Z}^n}
\langle k \rangle^{3/2} \left\|
\partial^\alpha_{x_\ell} \Box_k u \right\|_{L^\infty_t L^2_x  \bigcap  L^3_t L^{6}_{x}
(\mathbb{R}^{1+n})} \label{modX.m2}
\end{align}
and $\|u\|_{Y} \lesssim \delta$. Moreover, the scattering operator
of \eqref{gNLS} $S$ carries a zero neighborhood in $C(\mathbb{R},
M^{5/2}_{2,1})$ into $C(\mathbb{R}, M^{5/2}_{2,1})$.
\end{thm}

When the nonlinearity $F$ has a simple form, say,
\begin{align}
{\rm i} u_t + \Delta_\pm u = \sum^n_{i=1} \lambda_i
\partial_{x_i} (u^{\kappa_i+1} ), \quad u(0,x)= u_0(x),
\label{DNLS1}
\end{align}
we obtained in \cite{WaHu} the global well posedness of the DNLS
\eqref{DNLS1} for the small data in modulation spaces
$M^{1/\kappa_1}_{2,1}$ in one spatial dimension. In higher spatial
dimensions $n\ge 2$, we have

\begin{thm} \label{DNLS1-mod}
Let $n\ge 2$, $\kappa_i > 2, \; \kappa_i> 4/n$, $\kappa_i \in
\mathbb{N}$, $\lambda_i \in \mathbb{C}, $ $\kappa=\min_{1\le i \le
n} \, \kappa_i$. Assume that $ u_0 \in M^{1/2}_{2,1}$ and
$\|u_0\|_{M^{1/2}_{2,1}} \le \delta$ for some small $\delta>0$. Then
\eqref{DNLS1} has a unique global solution $u\in C(\mathbb{R},
M^{1/2}_{2,1}) \cap X$, where
\begin{align}
\|u\|_{X_1} =  & \sum^n_{i=1} \sum_{k\in \mathbb{Z}^n, \ |k_i|>4}
\langle k_i\rangle  \left\| \Box_k u \right\|_{L^\infty_{x_i}
L^2_{(x_j)_{j\not=i}}L^2_t(\mathbb{R}^{1+n})} \nonumber\\
&  + \sum^n_{i=1} \sum_{k\in \mathbb{Z}^n} \langle k
\rangle^{1/2-1/\kappa} \left\| \Box_k u \right\|_{L^{\kappa}_{x_i}
L^\infty_{(x_j)_{j\not=i}}L^\infty_t(\mathbb{R}^{1+n})} \nonumber\\
&  +  \sum_{k\in \mathbb{Z}^n}  \langle k \rangle^{1/2} \left\|
\Box_k u \right\|_{L^\infty_t L^2_x  \bigcap  L^{2+\kappa}_{x,t}
(\mathbb{R}^{1+n})} \label{modX.1}
\end{align}
and $\|u\|_{X_1} \lesssim \delta$.  Moreover, the scattering
operator of \eqref{DNLS1} $S$ carries a zero neighborhood in
$C(\mathbb{R}, M^{1/2}_{2,1})$ into $C(\mathbb{R}, M^{1/2}_{2,1})$.
\end{thm}
We remark that in Theorem \ref{DNLS1-mod}, the same result holds if
the nonlinear term $\partial_{x_i} (u^{\kappa_i+1} )$ is replaced by
$\partial_{x_i} (|u|^{\kappa_i}u )$ $(\kappa_i\in 2\mathbb{N})$.

\begin{thm} \label{DNLS1-modm2}
Let $n\ge 3$, $\kappa_i \in \mathbb{N}$, $\lambda_i \in \mathbb{C},
$ $\kappa=\min_{1\le i \le n} \, \kappa_i=2$. Assume that $ u_0 \in
M^{3/2}_{2,1}$ and $\|u_0\|_{M^{3/2}_{2,1}} \le \delta$ for some
small $\delta>0$. Then \eqref{DNLS1} has a unique global solution
$u\in C(\mathbb{R}, M^{3/2}_{2,1}) \cap Y_1$, where
\begin{align}
\|u\|_{Y_1} =  & \sum^n_{i=1} \sum_{k\in \mathbb{Z}^n, \ |k_i|>4}
\langle k_i\rangle^{2}  \left\| \Box_k u \right\|_{L^\infty_{x_i}
L^2_{(x_j)_{j\not=i}}L^2_t(\mathbb{R}^{1+n})} \nonumber\\
&  + \sum^n_{i=1} \sum_{k\in \mathbb{Z}^n}  \left\| \Box_k u
\right\|_{L^{2}_{x_i}
L^\infty_{(x_j)_{j\not=i}}L^\infty_t(\mathbb{R}^{1+n})} \nonumber\\
&  +  \sum_{k\in \mathbb{Z}^n}  \langle k \rangle^{3/2} \left\|
\Box_k u \right\|_{L^\infty_t L^2_x  \bigcap  L^{3}_{t}L^6_x
(\mathbb{R}^{1+n})}. \label{modY.1}
\end{align}
Moreover,  the scattering operator of \eqref{gNLS} $S$ carries a
zero neighborhood in $C(\mathbb{R}, M^{3/2}_{2,1})$ into
$C(\mathbb{R}, M^{3/2}_{2,1})$.
\end{thm}

\begin{cor} \label{DNLS-cor}
Let $n\ge 2$, $m\ge 2$, $s>(n+3)/2$. Assume that $ u_0 \in H^{s+1}$
and $\|u_0\|_{H^{s+1}} \le \delta$.  Then \eqref{gNLS} has a unique
global solution $u\in X$.
\end{cor}

When $m=1$, Christ \cite{Ch} showed the ill posedness of
\eqref{DNLS1} in any $H^s$ for one spatial dimension case. For
general nonlinearity in \eqref{gNLS}, we do not know what happens in
the case $m=1$ in higher spatial dimensions.

\subsection{Notations} \label{Notation}

The following are some notations which will be frequently used in
this paper: $\mathbb{C}, \mathbb{R}, \mathbb{N}$ and $ \mathbb{Z}$
will stand for the sets of complex number, reals, positive integers
and integers, respectively.   $c \le 1$, $C>1$ will denote positive
universal constants, which can be different at different places.
$a\lesssim b$ stands for $a\le C b$ for some constant $C>1$, $a\sim
b$ means that $a\lesssim b$ and $b\lesssim a$. We write $a\wedge b
=\min(a,b)$, $a\vee b =\max(a,b)$. We denote by $p'$ the dual number
of $p \in [1,\infty]$, i.e., $1/p+1/p'=1$. We will use Lebesgue
spaces $L^p:=L^p(\mathbb{R}^n)$, $\|\cdot\|_p :=\|\cdot\|_{L^p}$,
Sobolev spaces $H^{s}=(I-\Delta)^{-s/2}L^2$. Some properties of
these function spaces can be found in \cite{BL,Tr}. If there is no
explanation, we always assume that spatial dimensions $n\ge 2$. We
will use the function spaces $L^q_t L^p_x (\mathbb{R}^{n+1})$ and
$L^p_x L^q_t ( \mathbb{R}^{n+1})$ for which the norms are defined by
\begin{align}
& \|f\|_{L^q_{t} L^p_x (\mathbb{R}^{n+1})}= \left\|
\|f\|_{L^{p}_{x}(\mathbb{R}^{n})} \right\|_{L^{q}_t(\mathbb{R})} , \
\ \|f\|_{L^p_x L^q_{t}(\mathbb{R}^{n+1})}= \left\| \|f\|_{
L^{q}_t(\mathbb{R})} \right\|_{L^{p}_{x}(\mathbb{R}^{n})}, \nonumber
\end{align}
$L^p_{x,t}(\mathbb{R}^{n+1}):=L^p_x L^p_{t}(\mathbb{R}^{n+1})$. We
denote by $L^{p_1}_{x_i} L^{p_2}_{(x_j)_{j\not=i}} L^{p_2}_t:=
L^{p_1}_{x_i} L^{p_2}_{(x_j)_{j\not=i}} L^{p_2}_t(\mathbb{R}^{1+n})$
the anisotropic Lebesgue space for which the norm is defined by
\begin{align}
\|f\|_{L^{p_1}_{x_i} L^{p_2}_{(x_j)_{j\not=i}}L^{p_2}_t} =\left\|
\|f\|_{L^{p_2}_{x_1,...,x_{j-1}, x_{j+1},...,x_n}
L^{p_2}_t(\mathbb{R} \times \mathbb{R}^{n-1})}
\right\|_{L^{p_1}_{x_i}(\mathbb{R})} . \label{Notation.1}
\end{align}
It is also convenient to use the notation $L^{p_1}_{x_1}
L^{p_2}_{x_2,...,x_n}L^{p_2}_t: =L^{p_1}_{x_1}
L^{p_2}_{(x_j)_{j\not=1}}L^{p_2}_t $. For any $1<k<n$, we denote by
$\mathscr{F}_{x_1,...,x_k}$ the partial Fourier transform:
\begin{align}
(\mathscr{F}_{x_1,...,x_k} f)(\xi_1,...,\xi_k, x_{k+1},...,x_n) =
\int_{\mathbb{R}^k} e^{-{\rm
i}(x_1\xi_1+...+x_k\xi_k)}f(x)dx_1...dx_k \label{Notation.2}
\end{align}
and by $\mathscr{F}^{-1}_{\xi_1,...,\xi_k}$ the partial inverse
Fourier transform, similarly for $\mathscr{F}_{t,x}$ and
$\mathscr{F}^{-1}_{\tau,\xi}$. $\mathscr{F}:=
\mathscr{F}_{x_1,...,x_n}, \ \mathscr{F}^{-1}:=
\mathscr{F}^{-1}_{\xi_1,...,\xi_n} $. $D^s_{x_i}
=(-\partial^2_{x_i})^{s/2}= \mathscr{F}^{-1}_{\xi_i} |\xi_i|^s
\mathscr{F}_{x_i}$ expresses the partial Riesz potential in the
$x_i$ direction. $\partial^{-1}_{x_i} = \mathscr{F}^{-1}_{\xi_i}
({\rm i} \xi_i)^{-1} \mathscr{F}_{x_i}$. We will use the Bernstein
multiplier estimate; cf. \cite{BL,Tr}. For any $r\in [1,\infty]$,
\begin{align} \label{Bernstein}
\|\mathscr{F}^{-1}\varphi \mathscr{F} f\|_r\le
C\|\varphi\|_{H^s}\|f\|_r, \quad s>n/2.
\end{align}
We will use the frequency-uniform decomposition operators (cf.
\cite{Wa1,WaHe,WaHu}). Let $\{\sigma_k\}_{k\in \mathbb{Z}^n}$ be a
function sequence satisfying
\begin{align}
 \left\{\begin{array}{l}
\sigma_k(\xi)  \ge c, \quad \forall \; \xi \in Q_k,\\
 {\rm supp}\, \sigma_k \subset \{\xi: |\xi-k| \le \sqrt{n}\},\\
\sum_{k\in \mathbb{Z}^n} \sigma_k (\xi) \equiv 1, \quad \forall \;
\xi \in
\mathbb{R}^n,\\
|D^\alpha \sigma_k(\xi)| \le C_m, \quad \forall \; \xi \in \mathbb{
R}^n,\; |\alpha| \le m \in \mathbb{N}.
\end{array}\right. \label{UD-1}
\end{align}
Denote
\begin{align}
\Upsilon = \left\{ \{\sigma_k\}_{k\in \mathbb{Z}^n}: \;
\{\sigma_k\}_{k\in \mathbb{Z}^n}\;\; \mbox {satisfies \eqref{UD-1}}
\right\}. \label{funct-2}
\end{align}
Let $\{\sigma_k\}_{k\in \mathbb{Z}^n}\in \Upsilon$ be a function
sequence and
\begin{align}
\Box_k := \mathscr{F}^{-1} \sigma_k \mathscr{F}, \quad k\in \mathbb{
Z}^n, \label{funct-3}
\end{align}
which are said to be the frequency-uniform decomposition operators.
One may ask the existence of the frequency-uniform decomposition
operators. Indeed, let $\rho\in \mathscr{S}(\mathbb{R}^n)$ and
$\rho:\, \mathbb{R}^n\to [0,1]$ be a smooth radial bump function
adapted to the ball $B(0, \sqrt{n})$, say $\rho(\xi)=1$ as $|\xi|\le
\sqrt{n}/2$, and $\rho(\xi)=0$ as $|\xi| \ge \sqrt{n} $. Let
$\rho_k$ be a translation of $\rho$: $ \rho_k (\xi) = \rho (\xi- k),
\; k\in \mathbb{Z}^n$.  We write
\begin{align}
\eta_k (\xi)= \rho_k(\xi) \left(\sum_{k\in \mathbb{
Z}^n}\rho_k(\xi)\right)^{-1}, \quad k\in \mathbb{Z}^n.
\label{funct-1}
\end{align}
We have $\{\eta_k\}_{k\in \mathbb{Z}^n}\in \Upsilon$. It is easy to
see that for any $\{\eta_k\}_{k\in \mathbb{Z}^n}\in \Upsilon$,
$$
\|f\|_{M^s_{2,1}} \sim \sum_{k\in \mathbb{Z}^n} \langle k\rangle^s
\|\Box_k f\|_{L^2(\mathbb{R}^n)}.
$$
We will use the function space $\ell^{1,s}_\Box (L^p_{t}L^r_x
(I\times \mathbb{R}^n))$ which contains all of the functions
$f(t,x)$ so that the following norm is finite:
\begin{align} \|f\|_{\ell^{1,s}_\Box (L^p_{t}
L^r_x(I\times \mathbb{R}^n))}:= \sum_{k\in \mathbb{Z}^n} \langle
k\rangle ^s  \|\Box_k f\|_{L^p_{t}L^r_x (I\times \mathbb{R}^n)}.
\label{Mod.7}
\end{align}
For simplicity, we write  $\ell^{1}_\Box (L^p_{t}L^r_x (I\times
\mathbb{R}^n))= \ell^{1,0}_\Box (L^p_{t}L^r_x (I\times
\mathbb{R}^n))$.

This paper is organized as follows. In Section \ref{AGSE} we show
the global smooth effect estimates of the solutions of the linear
Schr\"odinger equation in anisotripic Lebesgue spaces. In Sections
\ref{LE1} and \ref{LE2} we consider the frequency-uniform localized
versions for the global maximal function estimates, the global
smooth effects, together with their relations to the Strichartz
estimates. In Sections \ref{pf-thm1} and \ref{pf-thm2} we prove our
Theorems \ref{DNLS1-mod} and \ref{DNLS-mod}, respectively. In the
Appendix we generalize the Christ-Kiselev Lemma to the anisotropic
Lebesgue spaces in higher dimensions.

\section{Anisotropic global smooth effects} \label{AGSE}

In this section, we always denote
$$
S(t)= e^{{\rm i}t \Delta_\pm}= \mathscr{F}^{-1} e^{{\rm i}t
\sum^n_{j=1}\varepsilon_j \xi^2_j}\mathscr{F}, \ \ \mathscr{A} f
(t,x)= \int^t_0 S(t-\tau) f(\tau,x) d\tau.
$$

\begin{prop} \label{GSE1}
For any $i=1,...,n$, we have the following estimate:
\begin{align}
\left\|\partial_{x_i} \mathscr{A}  f \right\|_{L^\infty_{x_i}
L^2_{(x_j)_{j\not= i}}L^2_t(\mathbb{R}^{1+n})} \lesssim
\|f\|_{L^1_{x_i} L^2_{(x_j)_{j\not= i}}L^2_t(\mathbb{R}^{1+n})}.
\label{smo-eff.1}
\end{align}
\end{prop}
{\bf Proof.} We have
\begin{align}
\partial_{x_1} \mathscr{A}  f  = c \mathscr{F}^{-1}_{t,x} \frac{\xi_1}{|\xi|^2_\pm
-\tau} \mathscr{F}_{t,x} f.  \label{smo-eff.2}
\end{align}
We can assume, without loss of generality that $|\xi|^2_\pm =
\xi^2_1+ \varepsilon_2 \xi^2_2+...+ \varepsilon_n \xi^n_n: =\xi^2_1+
|\bar{\xi}|^2_\pm$. By Plancherel's identity,
\begin{align}
& \left\|\partial_{x_1} \mathscr{A}  f \right\|_{L^\infty_{x_1}
L^2_{x_2,...,x_n} L^2_t(\mathbb{R}^{1+n})} \nonumber\\
& =\left\|\mathscr{F}^{-1}_{\xi_1} \frac{\xi_1}{\xi^2_1+
|\bar{\xi}|^2_\pm -\tau} \mathscr{F}_{t,x} f
\right\|_{L^\infty_{x_1} L^2_{\xi_2,...,\xi_n}
L^2_\tau(\mathbb{R}^{1+n})}  \nonumber\\
& \le \left\|\mathscr{F}^{-1}_{\xi_1} \frac{\xi_1}{\xi^2_1+
|\bar{\xi}|^2_\pm -\tau} \mathscr{F}_{t,x} f \right\|_{
L^2_{\xi_2,...,\xi_n}L^\infty_{x_1} L^2_\tau(\mathbb{R}^{1+n})}.
\label{smo-eff.3}
\end{align}
By changing the variable $\tau \to \mu+|\bar{\xi}|^2_\pm$, we have
\begin{align}
& \left\|\mathscr{F}^{-1}_{\xi_1} \frac{\xi_1}{\xi^2_1+
|\bar{\xi}|^2_\pm -\tau} \mathscr{F}_{t,x} f \right\|_{
L^2_{\xi_2,...,\xi_n}L^\infty_{x_1}
L^2_\tau(\mathbb{R}^{1+n})} \nonumber\\
& =  \left\|\mathscr{F}^{-1}_{\xi_1} \frac{\xi_1}{\xi^2_1-\mu}
\mathscr{F}_{t,x_1} (e^{-{\rm i}t |\bar{\xi}|^2_\pm}
\mathscr{F}_{x_2,...,x_n}f) \right\|_{
L^2_{\xi_2,...,\xi_n}L^\infty_{x_1} L^2_\mu(\mathbb{R}^{1+n})} .
\label{smo-eff.4}
\end{align}
Recalling the smooth effect estimate in one spatial dimension (cf.
\cite{KePoVe})
\begin{align}
\left\|\mathscr{F}^{-1}_{\tau, \xi} \frac{\xi}{\xi^2-\tau}
\mathscr{F}_{t,x} f \right\|_{L^\infty_{x} L^2_t(\mathbb{R}^{1+1})}
\lesssim \|f\|_{L^1_{x} L^2_t(\mathbb{R}^{1+1})}, \label{smo-eff.5}
\end{align}
we have from \eqref{smo-eff.3},  \eqref{smo-eff.4} and
\eqref{smo-eff.5} that
\begin{align}
& \left\|\partial_{x_1} \mathscr{A}  f \right\|_{L^\infty_{x_1}
L^2_{x_2,...,x_n} L^2_t(\mathbb{R}^{1+n})}
 \lesssim  \left\| e^{-{\rm i}t |\bar{\xi}|^2_\pm}
\mathscr{F}_{x_2,...,x_n}f \right\|_{ L^2_{\xi_2,...,\xi_n}L^1_{x_1}
L^2_t(\mathbb{R}^{1+n})}.  \label{smo-eff.6}
\end{align}
Using Minkowski's inequality and Plancherel's equality, we
immediately have
\begin{align}
& \left\|\partial_{x_1} \mathscr{A}  f \right\|_{L^\infty_{x_1}
L^2_{x_2,...,x_n} L^2_t(\mathbb{R}^{1+n})}
 \lesssim  \left\| f \right\|_{L^1_{x_1} L^2_{\xi_2,...,\xi_n}
L^2_t(\mathbb{R}^{1+n})}.  \label{smo-eff.7}
\end{align}
The other cases can be shown in a similar way. $\hfill\Box$

\begin{prop} \label{GSE2}
For any $i=1,...,n$, we have the following estimate:
\begin{align}
\left\| D^{1/2}_{x_i} S(t)u_0   \right\|_{L^\infty_{x_i}
L^2_{(x_j)_{j\not= i}}L^2_t(\mathbb{R}^{1+n})} \lesssim \|u_0\|_{2}.
\label{smo-eff.21}
\end{align}
\end{prop}
{\bf Proof.} By Plancherel's equality and Minkowski's inequality,
\begin{align}
\left\| S(t)u_0   \right\|_{L^\infty_{x_1}
L^2_{x_2,...,x_n}L^2_t(\mathbb{R}^{1+n})} &
=\left\|\mathscr{F}^{-1}_{\xi_1} e^{{\rm i}t\varepsilon_1 \xi^2_1}
\mathscr{F}_{x_1}(\mathscr{F}_{x_2,...,x_n} u_0)
\right\|_{L^\infty_{x_1}
L^2_{\xi_2,...,\xi_n}L^2_t(\mathbb{R}^{1+n})} \nonumber\\
& \le \left\|\mathscr{F}^{-1}_{\xi_1} e^{{\rm i}t\varepsilon_1
\xi^2_1} \mathscr{F}_{x_1}(\mathscr{F}_{x_2,...,x_n} u_0) \right\|_{
L^2_{\xi_2,...,\xi_n}L^\infty_{x_1}L^2_t(\mathbb{R}^{1+n})} .
\label{smo-eff.22}
\end{align}
Recall the half-order smooth effect of $S(t)$ in one spatial
dimension (cf. \cite{KePoVe}),
\begin{align}
\left\|\mathscr{F}^{-1}_{\xi} e^{{\rm i}t \xi^2} \mathscr{F}_{x} u_0
\right\|_{L^\infty_{x}L^2_t(\mathbb{R}^{1+1})} \lesssim
\|D^{-1/2}_{x} u_0\|_{L^2(\mathbb{R})}. \label{smo-eff.23}
\end{align}
Hence, in view of \eqref{smo-eff.22} and \eqref{smo-eff.23}, using
Plancherel's equality, we immediately have
\begin{align}
\left\| S(t)u_0   \right\|_{L^\infty_{x_1}
L^2_{x_2,...,x_n}L^2_t(\mathbb{R}^{1+n})} \lesssim \|D^{-1/2}_{x_1}
u_0\|_{L^2(\mathbb{R}^n)}, \label{smo-eff.24}
\end{align}
which implies the result, as desired. $\hfill \Box$

The dual version of \eqref{smo-eff.21} is
\begin{prop} \label{GSE3}
For any $i=1,...,n$, we have the following estimate:
\begin{align}
\left\| \partial_{x_i} \mathscr{A} f   \right\|_{L^\infty_{t}
L^2_{x} (\mathbb{R}^{1+n})} \lesssim \|D^{1/2}_{x_i} f\|_{L^1_{x_i}
L^2_{(x_j)_{j\not= i}}L^2_t(\mathbb{R}^{1+n})}. \label{smo-eff.31}
\end{align}
\end{prop}
{\bf Proof.} Denote $\mathbb{R}_+ =[0, \infty)$.  By Proposition
\ref{GSE2},
\begin{align}
& \left | \int_{\mathbb{R}_+} \left((\mathscr{A} \partial_{x_1} f)
(t) , \ \psi(t)
\right)dt \right| \nonumber\\
& = \left | \int_{\mathbb{R}_+} \left( f(\tau) , \ \int^\infty_\tau
S(\tau-t) \partial_{x_1}\psi(t) dt \right) d\tau \right| \nonumber\\
& \le \|  D^{1/2}_{x_1} f\|_{L^1_{x_1}
L^2_{x_2,...,x_n}L^2_t(\mathbb{R}^{1+n})} \int_{\mathbb{R}_+}
\|\partial_{x_1} S(\tau-t)D^{-1/2}_{x_1} \psi(t)\|_{L^\infty_{x_1}
L^2_{x_2,...,x_n}L^2_\tau(\mathbb{R}^{1+n})} dt \nonumber\\
& \lesssim  \| D^{1/2}_{x_1} f\|_{L^1_{x_1}
L^2_{x_2,...,x_n}L^2_t(\mathbb{R}^{1+n})} \|   \psi \|_{L^1_{t}
L^2_{x}(\mathbb{R}^{1+n})}.
 \label{smo-eff.32}
\end{align}
By duality, we have the result. $\hfill\Box$

\section{Linear estimates with $\Box_k$-decomposition} \label{LE1}

In this section we consider the smooth effect estimates, the maximal
function estimates, the Strichartz estimates and their interaction
estimates for the solutions of the linear Schr\"odinger equations by
using the frequency-uniform decomposition operators. For
convenience, we will use the following function sequence
$\{\sigma_k\}_{k\in \mathbb{Z}^n}$:

\begin{lem}\label{eq-mod-funct}
Let $\eta_k : \mathbb{R} \to [0,1]$ $(k\in \mathbb{Z})$ be a
smooth-function sequence satisfying condition \eqref{UD-1}. Denote
\begin{align}
\sigma_k (\xi):= \eta_{k_1}(\xi_1)...\eta_{k_n}(\xi_n), \ \
k=(k_1,...,k_n). \label{UD-2}
\end{align}
Then we have $\{\sigma_k\}_{k\in \mathbb{Z}^n} \in \Upsilon$.
\end{lem}
Recall that in \cite{WaHe}, we established the following Strichartz
estimates in a class of function spaces by using the
frequency-uniform decomposition operators.
\begin{lem}\label{Strichartz-mod}
 Let $2\le p< \infty$, $\gamma \ge 2\vee \gamma(p)$,
\begin{align}
 \frac{2}{\gamma(p)}=n\Big(\frac{1}{2}-\frac{1}{p}\Big). \nonumber
\end{align}
 Then we have
\begin{align}
\left\|S(t) \varphi \right\|_{\ell^1_\Box (L^\gamma(\mathbb{R},
L^p(\mathbb{R}^n))
)} & \lesssim \| \varphi\|_{M_{2, 1}(\mathbb{R}^n)},  \nonumber\\
\left\|\mathscr{A} f \right\| _{\ell^1_\Box (L^\gamma (\mathbb{R},
L^p(\mathbb{R}^n))) \cap \ell^1_\Box (L^\infty(\mathbb{R},
L^2(\mathbb{R}^n))) } & \lesssim \|f\|_{\ell^1_\Box
(L^{\gamma'}(\mathbb{R}, L^{p'}(\mathbb{R}^n)))}. \nonumber
\end{align}
In particular, if $2+4/n\le p< \infty$, then we have
\begin{align}
\left\|S(t) \varphi \right\|_{\ell^1_\Box
(L^p_{t,x}(\mathbb{R}^{1+n})
)} & \lesssim \| \varphi\|_{M_{2, 1}(\mathbb{R}^n)},  \nonumber\\
\left\|\mathscr{A} f \right\| _{\ell^1_\Box
(L^p_{t,x}(\mathbb{R}^{1+n})) \, \cap \, \ell^1_\Box
(L^\infty_tL^2_x(\mathbb{R}^{1+n})) } & \lesssim \|f\|_{\ell^1_\Box
(L^{p'}_{t,x}(\mathbb{R}^{1+n}))}. \nonumber
\end{align}
\end{lem}

The next lemma is essentially known, see \cite{Tr,Wa1}.

\begin{lem}\label{qnorm:pnorm}
Let $\Omega\subset {\Bbb R}^n$ be a compact set with ${\rm
diam}\,\Omega <2R$, $0<p\le q \le \infty.$ Then there exists a
constant $C>0$, which depends only on $p,q$ such that
\begin{align}
\| f\|_q \le C R^{n(1/p-1/q)}\|f\|_p, \quad \forall\; f\in
L^p_\Omega, \nonumber
\end{align}
where $L^p_\Omega =\{f\in {\cal S}'({\Bbb R}^n): {\rm supp}
\hat{f}\subset \Omega, \; \|f\|_p<\infty\}.$
\end{lem}

In Lemma \ref{qnorm:pnorm} we emphasize that the constant $C>0$  is
independent of the position of $\Omega$ in frequency spaces, say, in
the case $\Omega=B(k, \sqrt{n}),  \ k\in {\Bbb Z}^n$, Lemma
\ref{qnorm:pnorm} uniformly holds for all $k\in {\Bbb Z}^n$.

\begin{lem}\label{discret-deriv}
We have for any $\sigma \in  \mathbb{R}$ and  $k=(k_1,...,k_n)\in
\mathbb{Z}^n$ with $|k_i| \ge 4$,
\begin{align*}
\|\Box_k D^\sigma_{x_i} u\|_{L^{p_1}_{x_1} L^{p_2}_{x_2,...,x_n}
L^{p_2}_t (\mathbb{R}^{1+n}) } & \lesssim \langle k_i
\rangle^{\sigma} \|\Box_k u\|_{L^{p_1}_{x_1} L^{p_2}_{x_2,...,x_n}
L^{p_2}_t (\mathbb{R}^{1+n}) }.
\end{align*}
Replacing $D^\sigma_{x_i}$ by $\partial^\sigma_{x_i}$ ($\sigma\in
\mathbb{N}$), the above inequality holds for all $k\in
\mathbb{Z}^n$.
\end{lem}
{\bf Proof.} Using Lemma \ref{eq-mod-funct}, one has that
\begin{align*}
\Box_k D^\sigma_{x_i} u = \sum^{1}_{\ell=-1} \int_{\mathbb{R}}
\left(\mathscr{F}^{-1}_{\xi_i} (\eta_{k_i+\ell}(\xi_i)
|\xi_i|^\sigma) \right)(y_i)  (\Box_k u)(x_i-y_i) dy_i.
\end{align*}
It follows that
\begin{align*}
& \|\Box_k D^\sigma_{x_i} u\|_{L^{p_1}_{x_1} L^{p_2}_{x_2,...,x_n}
L^{p_2}_t (\mathbb{R}^{1+n}) } \\
& \lesssim \sum^{1}_{\ell=-1}  \|\mathscr{F}^{-1}_{\xi_i}
(\eta_{k_i+\ell}(\xi_i) |\xi_i|^\sigma)\|_{L^1(\mathbb{R})} \|
\Box_k u\|_{L^{p_1}_{x_1}
L^{p_2}_{x_2,...,x_n} L^{p_2}_t (\mathbb{R}^{1+n}) }\\
& \lesssim  \langle k_i \rangle^{\sigma} \| \Box_k
u\|_{L^{p_1}_{x_1} L^{p_2}_{x_2,...,x_n} L^{p_2}_t
(\mathbb{R}^{1+n}) }.
\end{align*}
The result follows. $\hfill \Box$

\medskip

Ionescu and Kenig \cite{Io-Ke} showed the following maximal function
estimates in higher spatial dimensions $n\ge 3$:
\begin{align}
 \|\triangle_k S(t)
u_0\|_{L^{2}_{x_i} L^\infty_{(x_j)_{j\not=i}} L^\infty_t
(\mathbb{R}^{1+n}) } & \lesssim  2^{(n-1)k/2} \| \triangle_k
u_0\|_{L^2(\mathbb{R}^n)}. \label{MFM-2a}
\end{align}
We partially resort to their idea to obtain the following

\begin{prop}\label{MaxFunct-Mod}
 Let $4/n< q \le \infty$, $q\ge 2$. Then we have
\begin{align}
 \|\Box_k S(t)
u_0\|_{L^{q}_{x_i} L^\infty_{(x_j)_{j\not=i}} L^\infty_t
(\mathbb{R}^{1+n}) } & \lesssim  \langle k_i \rangle^{1/q} \| \Box_k
u_0\|_{L^2(\mathbb{R}^n)}. \label{MFM-2}
\end{align}
\end{prop}
{\bf Proof.} For convenience, we write $\bar{x}=(x_1,...,x_{n-1})$.
By duality, it suffices to show that for any $\varphi\in
L^{q'}_{x_1}L^1_{\bar{x},t} (\mathbb{R}^{1+n}) \bigcap \mathscr{S}
(\mathbb{R}^{1+n})$ with $\varphi(t)=\pm \varphi(-t)$,
\begin{align}
 \int_{\mathbb{R}} (\Box_k S(t) u_0, \ \varphi(t))dt   \lesssim  \langle k_i \rangle^{1/q} \| \Box_k
u_0\|_{L^2(\mathbb{R}^n)}
\|S(t)\varphi\|_{L^{q'}_{x_1}L^1_{\bar{x},t}(\mathbb{R}^{1+n})}.
\label{MFM-a}
\end{align}
By duality, we have
\begin{align}
 \int_{\mathbb{R}} (\Box_k S(t) u_0, \ \varphi(t))dt   \lesssim   \|
u_0\|_{L^2(\mathbb{R}^n)} \left\|\int_{\mathbb{R}} \Box_k S(-t)
\varphi(t)dt \right\|_{L^{2}(\mathbb{R}^n)}. \label{MFM-b}
\end{align}
We have from Lemma \ref{discret-deriv} that
\begin{align}
\left\|\int_{\mathbb{R}} \Box_k S(-t) \varphi(t)dt
\right\|^2_{L^{2}(\mathbb{R}^n)} \lesssim
\|S(t)\varphi\|_{L^{q'}_{x_1}L^1_{\bar{x},t} (\mathbb{R}^{1+n})}
\left\|\int_{\mathbb{R}} \Box_k S(2t-\tau) \varphi(\tau)d\tau
\right\|_{L^{q}_{x_1} L^\infty_{\bar{x},t} (\mathbb{R}^{1+n})}.
\label{MFM-c}
\end{align}
In view of Lemma \ref{eq-mod-funct}, we can write $\Box_k=
\mathscr{F}^{-1} \eta_{k_1}(\xi_1)...\eta_{k_n}(\xi_n)\mathscr{F}:=
\mathscr{F}^{-1}
\eta_{k_1}(\xi_1)\eta_{\bar{k}}(\bar{\xi})\mathscr{F}$. By
Minkowski's and Young's inequalities,
\begin{align}
&  \left\|\int_{\mathbb{R}} \Box_k S(2t-\tau) \varphi(\tau)d\tau
\right\|_{L^{q}_{x_1}L^\infty_{\bar{x},t}(\mathbb{R}^{1+n})} \nonumber\\
& \lesssim \left \| \mathscr{F}^{-1} e^{{\rm i} 2t|\xi|^2_{\pm}}
\eta_{\bar{k}}(\bar{\xi}) \eta_{k_1}(\xi_1) \mathscr{F}
\int_{\mathbb{R}}S(-\tau) \varphi(\tau) d\tau \right\|_{L^{q}_{x_1}
L^\infty_{\bar{x},t} (\mathbb{R}^{n+1})}
\nonumber\\
& \lesssim   \left\| \|\mathscr{F}^{-1}_{\bar{\xi}} e^{{\rm i}
2t|\bar{\xi}|^2_{\pm}} \eta_{\bar{k}}(\bar{\xi})
\|_{L^\infty_{\bar{x}} (\mathbb{R}^{n})}
\left\|\mathscr{F}^{-1}_{\xi_1} e^{{\rm i} 2t\xi^2_1}
\eta_{k_1}(\xi_1) \mathscr{F}_{x_1}
 \int_{\mathbb{R}} S(-\tau)\varphi (\tau) d\tau \right\|_{L^1_{\bar{x}}}
\right\|_{L^{q}_{x_1}L^{\infty}_{t} }
\nonumber\\
& \lesssim  \|\mathscr{F}^{-1} e^{{\rm i} 2t|\xi|^2_{\pm}}
\eta_{k_1}(\xi_1)\eta_{\bar{k}}(\bar{\xi})
\|_{L^{q/2}_{x_1}L^{\infty}_{\bar{x},t} (\mathbb{R}^{n+1})}
\left\|\int_{\mathbb{R}} S(-\tau) \varphi (\tau)d\tau
\right\|_{L^{q'}_{x_1} L^1_{\bar{x}}(\mathbb{R}^{n})}
\nonumber\\
& \lesssim  \|\mathscr{F}^{-1} e^{{\rm i}2 t|\xi|^2_{\pm}}
\eta_{k_1}(\xi_1)\eta_{\bar{k}}(\bar{\xi})
\|_{L^{q/2}_{x_1}L^{\infty}_{\bar{x},t} (\mathbb{R}^{1+n})}
\|S(\tau)
 \varphi\|_{L^{q'}_{x_1}L^1_{\bar{x},t}(\mathbb{R}^{n+1})}.
 \label{MFM-d}
\end{align}
Hence, it suffices to show that
\begin{align}
\|\mathscr{F}^{-1} e^{{\rm i} t|\xi|^2_{\pm}}
\eta_{k_1}(\xi_1)\eta_{\bar{k}}(\bar{\xi})
\|_{L^{q/2}_{x_1}L^{\infty}_{\bar{x},t} (\mathbb{R}^n)} \lesssim
\langle k_1\rangle^{2/q}. \nonumber
\end{align}
In view of the decay of $\Box_k S(t)$, we see that (cf. \cite{WaHe})
\begin{align}
& \|\mathscr{F}^{-1}_{\bar{\xi}} e^{{\rm i} t|\bar{\xi}|^2_{\pm}}
\eta_{\bar{k}}(\bar{\xi}) \|_{L^\infty_{\bar{x}} (\mathbb{R}^{n-1})}
\lesssim (1+|t|)^{-(n-1)/2}, \nonumber\\
& \|\mathscr{F}^{-1}_{\xi_1} e^{{\rm i} t \xi_1^2} \eta_{k_1}(\xi_1)
\|_{L^\infty_{x_1} (\mathbb{R})} \lesssim (1+|t|)^{-1/2}. \nonumber
\end{align}
On the other hand, integrating by part, one has that for $|x_1| >
4|t|\langle k_1\rangle$,
\begin{align}
& |\mathscr{F}^{-1}_{\xi_1} e^{{\rm i} t \xi_1^2} \eta_{k_1}(\xi_1)
| \lesssim |x_1|^{-2}. \nonumber
\end{align}
Hence, for $|x_1|>1$,
\begin{align}
& |\mathscr{F}^{-1} e^{{\rm i} t |\xi|_{\pm}^2} \eta_{k_1}(\xi_1)
\eta_{\bar{k}}(\bar{\xi})| \lesssim (1+|x_1|)^{-2} + \langle
k_1\rangle^{n/2}(\langle k_1\rangle +|x_1|)^{-n/2}. \nonumber
\end{align}
So, we have
\begin{align}
 \|\mathscr{F}^{-1} e^{{\rm i} t |\xi|_{\pm}^2} \eta_{k_1}(\xi_1)
\eta_{\bar{k}}(\bar{\xi})\|_{L^{q/2}_{x_1}L^{\infty}_{\bar{x},t}
(\mathbb{R}^n)} & \lesssim 1 + \langle k_1\rangle^{n/2}\|(\langle
k_1\rangle +|x_1|)^{-n/2}\|_{L^{q/2}_{x_1} (\mathbb{R})} \nonumber\\
& \lesssim \langle k_1\rangle^{2/q}. \nonumber
\end{align}
This finishes the proof of  \eqref{MFM-2}. $\hfill \Box$

\begin{rem} \rm We conjecture that \eqref{MFM-2} also holds in the
case $p=4/n$ if $n=2$.  We now show that \eqref{MFM-2} is sharp.
Indeed, take $k_1\in \mathbb{N}$,  $\mathscr{F}_x u_0(\xi) =
\eta_{k_1}(\xi_1) \eta_0(\xi_2)...\eta_0(\xi_n)$, where $\eta_k$ is
as in Lemma \ref{eq-mod-funct}. For $\sigma_k:=
\eta_{k_1}(\xi_1)\eta_0(\xi_2)...\eta_0(\xi_n)$, we easily see that
\begin{align}
& \|\Box_{k} S(t)u_{0}\|^{q}_{L^q_{x_1}L^{\infty}_{x_2,...,x_n}
L^\infty_t (\mathbb{R}^{1+n})} \nonumber\\
& = \left\|\mathscr{F}^{-1}_{\xi_1} ( e^{{\rm i} \varepsilon_1 t
\xi^2_1} \eta^2_{k_1}(\xi_1)) \prod^n_{i=2} \mathscr{F}^{-1}_{\xi_i}
( e^{{\rm i} \varepsilon_i t \xi^2_i} \eta^2_{0}(\xi_i))
\right\|^{q}_{L^q_{x_1}L^{\infty}_{x_2,...,x_n}
L^\infty_t (\mathbb{R}^{1+n})} \nonumber\\
& \ge \int_{|x_1| \le c k_1}  \sup_{t,
x_2,...,x_n}\left|\mathscr{F}^{-1}_{\xi_1} ( e^{{\rm i}
\varepsilon_1 t \xi^2_1} \eta^2_{k_1}(\xi_1)) (x_1) \prod^n_{i=2}
\mathscr{F}^{-1}_{\xi_i} ( e^{{\rm i} \varepsilon_i t \xi^2_i}
\eta^2_{0}(\xi_i)) (x_i) \right|^{q} dx_1.  \nonumber
\end{align}
Taking $t= -x_1/2 \varepsilon_1 k_1$ and $|x_i|< c$, $i=2,...,n$, we
easily see that
\begin{align}
& \left|\mathscr{F}^{-1}_{\xi_1} ( e^{{\rm i} \varepsilon_1 t
\xi^2_1} \eta^2_{k_1}(\xi_1)) (x_1) \right| \gtrsim 1,   \nonumber\\
& \left| \mathscr{F}^{-1}_{\xi_i} ( e^{{\rm i} \varepsilon_i t
\xi^2_i} \eta^2_{0}(\xi_i)) (x_i) \right| \gtrsim 1.   \nonumber
\end{align}
Therefore, we have
\begin{align}
& \|\Box_{k} S(t)u_{0}\|^{q}_{L^q_{x_1}L^{\infty}_{x_2,...,x_n}
L^\infty_t (\mathbb{R}^{1+n})}  \gtrsim k_1. \nonumber
\end{align}
\end{rem}
The dual version of Proposition \ref{MaxFunct-Mod} is the following
\begin{prop}\label{MaxFunct-Mod-Dual}
 Let $2\le  q \le \infty$, $q>4/n$. Then we have for any $k=(k_1,...,k_n)\in
 \mathbb{Z}^n$,
\begin{align}
\left \|\Box_k \int_\mathbb{R} S(t-\tau) f(\tau) d\tau
\right\|_{L^\infty_t L^2_x(\mathbb{R}^{1+n})} & \lesssim \langle k_i
\rangle^{1/q} \| \Box_k f\|_{L^{q'}_{x_i} L^1_{(x_j)_{j\not=i}}
L^1_t (\mathbb{R}^{1+n}) }. \label{MFM-2-dual}
\end{align}
\end{prop}
{\bf Proof.} Denote
\begin{align}
\tilde{\Box}_k = \sum_{\ell \in \Lambda} \Box_{k+\ell}, \ \ \Lambda=
\{\ell\in \mathbb{Z}^n: \ {\rm supp}\, \sigma_k \cap {\rm supp}\,
\sigma_{k+\ell} \not=\varnothing \}.  \label{ma-sm-mo-3}
\end{align}
Write
\begin{align}
 \mathcal {L}_k (f, \psi) &:= \left|\int_\mathbb{R}\left( \Box_k
\int_{\mathbb{R}}S(t-\tau)  f(\tau) d\tau, \  \psi(t) \right ) dt
\right|
 \label{ma-sm-mo-n}
\end{align}
By Proposition \ref{MaxFunct-Mod},
\begin{align}
 \mathcal {L}_k (f, \psi) & = \left|\left( \Box_k
f(\tau) , \  \tilde{\Box}_k \int_\mathbb{R} S(\tau-t) \psi(t) dt \right ) d\tau \right| \nonumber\\
& \le \left\| \Box_k f \right\|_{L^{q'}_{x_i} L^1_{(x_j)_{j\not=i}}
L^1_t (\mathbb{R}^{1+n}) } \left\|\tilde{\Box}_k \int_\mathbb{R}
S(\tau -t) \psi(t) dt \right\|_{L^{q}_{x_i}
L^\infty_{(x_j)_{j\not=i}} L^\infty_t (\mathbb{R}^{1+n})}
\nonumber\\
& \le \left\| \Box_k f \right\|_{L^{q'}_{x_i} L^1_{(x_j)_{j\not=i}}
L^1_t (\mathbb{R}^{1+n}) } \langle k_i \rangle^{1/q}
\left\|\tilde{\Box}_k \psi \right\|_{L^{1}_{t} L^2_{x}
(\mathbb{R}^{1+n})}.
 \label{ma-sm-mo-4}
\end{align}
By duality, we have the result, as desired. $\hfill \Box$

In view of Propositions \ref{GSE1} and \ref{GSE3}, we have
\begin{prop}\label{1order-sm-mod}
We have for any $k=(k_1,...,k_n)\in
 \mathbb{Z}^n$,
\begin{align}
& \left \|\Box_k \mathscr{A} \partial_{x_i} f
\right\|_{L^\infty_{x_i} L^2_{(x_j)_{j\not=i}} L^2_t
(\mathbb{R}^{1+n})}  \lesssim  \|\Box_k f\|_{L^1_{x_i}
L^2_{(x_j)_{j\not=i}} L^2_t
(\mathbb{R}^{1+n}) },  \label{1-sm-mod} \\
& \left \|\Box_k \mathscr{A} \partial_{x_i} f \right\|_{L^\infty_t
L^2_x(\mathbb{R}^{1+n})} \lesssim  \langle k_i\rangle^{1/2} \|
\Box_k f\|_{L^1_{x_i} L^2_{(x_j)_{j\not=i}} L^2_t (\mathbb{R}^{1+n})
}. \label{1/2-sm-mod}
\end{align}
\end{prop}
{\bf Proof.} By Proposition \ref{GSE1}, we immediately have
\eqref{1-sm-mod}. In view of Proposition \ref{GSE3} and Lemma
\ref{discret-deriv}, we have \eqref{1/2-sm-mod} in the case $|k_i|
\ge 3$. If $|k_i| \le 2$, in view of Proposition \ref{GSE3},
\begin{align*}
& \left \|\Box_k \mathscr{A} \partial_{x_i} f \right\|_{L^\infty_t
L^2_x(\mathbb{R}^{1+n})} \lesssim  \left \|D^{-1/2}_{x_i} \Box_k
\mathscr{A} \partial_{x_i} f \right\|_{L^\infty_t
L^2_x(\mathbb{R}^{1+n})} \lesssim  \| \Box_k f\|_{L^1_{x_i}
L^2_{(x_j)_{j\not=i}} L^2_t (\mathbb{R}^{1+n}) },
\end{align*}
which implies the result, as desired. $\hfill \Box$

By the duality, we also have the following
\begin{prop}\label{ma-sm-mo}
 Let $2<  q \le \infty$ $q>4/n$. Then we have
\begin{align}
 \left \|\Box_k \mathscr{A} \partial_{x_i} f \right\|_{L^{q}_{x_i}
L^\infty_{(x_j)_{j\not=i}} L^\infty_t (\mathbb{R}^{1+n}) } &
\lesssim \langle k_i\rangle^{1/2+ 1/q}  \| \Box_k f\|_{L^1_{x_i}
L^2_{(x_j)_{j\not=i}} L^2_t (\mathbb{R}^{1+n}) }. \label{ma-sm-mo-2}
\end{align}
\end{prop}
{\bf Proof.} By Propositions \ref{MaxFunct-Mod-Dual},
\ref{1order-sm-mod} and Lemma \ref{discret-deriv},
\begin{align}
 \mathcal {L}_k (\partial_{x_1} f, \psi)
& = \left|\left( \Box_k \int_\mathbb{R} S(-\tau)
\partial_{x_1} f(\tau) d\tau, \  \tilde{\Box}_k \int_\mathbb{R} S(-t) \psi(t) dt \right )  \right| \nonumber\\
& \le \left\| \Box_k \int_\mathbb{R} S(-\tau)
\partial_{x_1} f(\tau) d\tau \right\|_{L^2(\mathbb{R}^n)}
\left\|  \tilde{\Box}_k \int_\mathbb{R} S(-t) \psi(t) dt   \right\|_{L^2(\mathbb{R}^n)} \nonumber\\
& \lesssim  \langle k_1\rangle^{1/2}  \|\Box_k  f\|_{L^{1}_{x_1}
L^2_{x_2,...,x_n}L^2_t(\mathbb{R}^{1+n})} \langle k_1 \rangle^{1/q}
\left\| \tilde{\Box}_k \psi \right\|_{L^{q'}_{x_1}
L^1_{x_2,...,x_n}L^1_t(\mathbb{R}^{1+n})}\nonumber\\
 & \lesssim \langle k_1\rangle^{1/2+1/q}\|\Box_k  f\|_{L^{1}_{x_1}
L^2_{x_2,...,x_n}L^2_t(\mathbb{R}^{1+n})} \|\psi \|_{L^{q'}_{x_1}
L^1_{x_2,...,x_n}L^1_t(\mathbb{R}^{1+n})}.
 \label{ma-sm-mo-5}
\end{align}
Again, by duality,  it follows from \eqref{ma-sm-mo-5} and
Christ-Kiselev's Lemma that \eqref{ma-sm-mo-2}  holds.
 $\hfill\Box$

\begin{prop}\label{st-sm-mo}
 Let $2\le r < \infty$, $2/\gamma(r)=n(1/2-1/r)$ and $\gamma> \gamma(r)\vee 2$.  We have
\begin{align}
& \left \|\Box_k S(t) u_0 \right\|_{ L^\gamma_t
L^r_x(\mathbb{R}^{1+n}) } \lesssim   \|\Box_k u_0\|_{L^2
(\mathbb{R}^{n}) },  \label{st-sm-mo-1} \\
& \left \|\Box_k \mathscr{A}  f \right\|_{L^\infty_t L^2_x \, \cap
\, L^\gamma_t L^r_x (\mathbb{R}^{1+n}) } \lesssim    \|\Box_k f\|_{
L^{\gamma'}_t L^{r'}_x (\mathbb{R}^{1+n})
}, \label{st-sm-mo-2}\\
& \left \|\Box_k \mathscr{A} \partial_{x_i}  f \right\|_{L^\gamma_t
L^r_x (\mathbb{R}^{1+n}) } \lesssim \langle k_i\rangle^{1/2} \|
\Box_k f\|_{L^1_{x_i} L^2_{(x_j)_{j\not=i}} L^2_t
(\mathbb{R}^{1+n})}, \label{st-sm-mo-3}\\
& \left \|\Box_k \mathscr{A} \partial_{x_i}  f
\right\|_{L^\infty_{x_i} L^2_{(x_j)_{j\not=i}} L^2_t
(\mathbb{R}^{1+n})}   \lesssim \langle k_i\rangle^{1/2}\| \Box_k
f\|_{ L^{\gamma'}_t L^{r'}_x (\mathbb{R}^{1+n})}, \label{st-sm-mo-4}
\end{align}
and for $2\le q<\infty$, $q>4/n$, $\alpha=0,1$,
\begin{align}
\left \|\Box_k \mathscr{A} \partial^\alpha_{x_i}  f
\right\|_{L^q_{x_i} L^\infty_{(x_j)_{j\not=i}} L^\infty_t
(\mathbb{R}^{1+n})} & \lesssim \langle k_i\rangle^{\alpha+ 1/q}
 \| \Box_k f\|_{ L^{\gamma'}_t L^{r'}_x
(\mathbb{R}^{1+n})}, \label{st-sm-mo-6}
\end{align}
\end{prop}
{\bf Proof.} From Lemma \ref{Strichartz-mod} it follows that
\eqref{st-sm-mo-1} and \eqref{st-sm-mo-2} hold. We now show
\eqref{st-sm-mo-3}. We use the same notations as in Proposition
\ref{ma-sm-mo}. By Lemmas \ref{Strichartz-mod}, \ref{discret-deriv}
and Proposition \ref{1order-sm-mod},
\begin{align}
 \mathcal {L}_k (\partial_{x_1} f, \psi) & \lesssim \langle k_i\rangle^{1/2} \|\Box_k  f\|_{L^{1}_{x_1}
L^2_{x_2,...,x_n}L^2_t(\mathbb{R}^{1+n})} \left\| \tilde{\Box}_k
\psi \right\|_{L^{\gamma'}_t L^{r'}_{x} (\mathbb{R}^{1+n})}
\nonumber\\
& \lesssim \langle k_i\rangle^{1/2} \|\Box_k f\|_{L^{1}_{x_1}
L^2_{x_2,...,x_n}L^2_t(\mathbb{R}^{1+n})} \left\| \psi
\right\|_{L^{\gamma'}_t L^{r'}_{x} (\mathbb{R}^{1+n})}.
 \label{st-sm-mo-7}
\end{align}
By duality,  it follows from \eqref{st-sm-mo-7} and Christ-Kiselev's
Lemma that \eqref{st-sm-mo-3}  holds. Exchanging the roles of $f$
and $\psi$, we immediately have \eqref{st-sm-mo-4} in the case
$r>2$. If $r=2$, \eqref{st-sm-mo-4} is a straightforward consequence
of the $1/2$-order smooth effect of $S(t)$.  By Lemmas
\ref{Strichartz-mod}, \ref{discret-deriv}, Proposition
\ref{MaxFunct-Mod-Dual}, and Christ-Kiselev's Lemma that we have
\eqref{st-sm-mo-6} in the case $q>2$, or $q=2$ and $r>2$. In the
case $q=r=2$, in view of the maximal function estimate, we see that
\eqref{st-sm-mo-6} also holds. $\hfill\Box$

\begin{cor}\label{st-sm-m-c}
 Let $4/n \le p < \infty$, $2\le q < \infty$, $q>4/n$.  We have
\begin{align}
& \left \| D^{1/2}_{x_1} \Box_k S(t) u_0 \right\|_{L^\infty_{x_1}
L^2_{x_2,...,x_n} L^2_t (\mathbb{R}^{1+n})}   \lesssim \| \Box_k
u_0\|_{ L^2 (\mathbb{R}^{n})}, \label{st-sm-m-c2}\\
& \left \| \Box_k S(t) u_0 \right\|_{L^q_{x_1}
L^\infty_{x_2,...,x_n} L^\infty_t (\mathbb{R}^{1+n})}   \lesssim
\langle k_i \rangle^{1/q} \| \Box_k u_0\|_{ L^2 (\mathbb{R}^{n})},
\label{st-sm-m-c3} \\
 & \left \|\Box_k S(t) u_0 \right\|_{
L^{2+p}_{t,x}\, \cap \, L^\infty_t L^2_x (\mathbb{R}^{1+n}) }
\lesssim \|\Box_k u_0\|_{L^2 (\mathbb{R}^{n}) }, \label{st-sm-m-c1}
\end{align}
\begin{align}
& \left \|\Box_k \mathscr{A} \partial_{x_1} f
\right\|_{L^\infty_{x_1} L^2_{x_2,...,x_n} L^2_t (\mathbb{R}^{1+n})}
\lesssim  \| \Box_k f\|_{L^1_{x_1} L^2_{x_2,...,x_n} L^2_t
(\mathbb{R}^{1+n})},
\label{st-sm-mo-5}\\
& \left \|\Box_k \mathscr{A} \partial_{x_1}  f \right\|_{L^q_{x_1}
L^\infty_{x_2,...,x_n} L^\infty_t (\mathbb{R}^{1+n})}  \lesssim
\langle k_1\rangle^{1/2+1/q} \| \Box_k f\|_{L^1_{x_1}
L^2_{x_2,...,x_n} L^2_t (\mathbb{R}^{1+n})},
\label{st-sm-m-c6}\\
& \left \|\Box_k \mathscr{A}  f \right\|_{L^\infty_t L^2_x \, \cap
\, L^{2+p}_{t,x}  (\mathbb{R}^{1+n}) } \lesssim  \langle
k_1\rangle^{1/2}  \|\Box_k f\|_{ L^1_{x_1} L^2_{x_2,...,x_n} L^2_t
(\mathbb{R}^{1+n}) }. \label{st-sm-m-c4}
\end{align}
\begin{align}
&\left \|\Box_k \mathscr{A} \partial_{x_1} f
\right\|_{L^\infty_{x_1} L^2_{x_2,...,x_n} L^2_t (\mathbb{R}^{1+n})}
\lesssim \langle k_1\rangle^{1/2}\| \Box_k f\|_{
L^{(2+p)/(1+p)}_{t,x} (\mathbb{R}^{1+n})},
\label{st-sm-mo-7}\\
& \left \|\Box_k \mathscr{A} \partial_{x_1}  f \right\|_{L^q_{x_1}
L^\infty_{x_2,...,x_n} L^\infty_t (\mathbb{R}^{1+n})}  \lesssim
 \langle k_1 \rangle^{1+1/q} \| \Box_k
f\|_{L^{(2+p)/(1+p)}_{t,x} (\mathbb{R}^{1+n})},
\label{st-sm-m-c8}\\
& \left \|\Box_k \mathscr{A}  f \right\|_{L^\infty_t L^2_x \, \cap
\, L^{2+p}_{t,x}  (\mathbb{R}^{1+n}) } \lesssim    \|\Box_k f\|_{
L^{(2+p)/(1+p)}_{t,x} (\mathbb{R}^{1+n}) }. \label{st-sm-m-c9}
\end{align}
In \eqref{st-sm-m-c6}, $q>2$ is required. Moreover, replacing
$L^{2+p}_{x,t}$ by $L^3_tL^6_x$, the results also hold.
\end{cor}

\section{Linear estimates with derivative interaction} \label{LE2}

In view of \eqref{st-sm-mo-5} in Corollary \ref{st-sm-m-c}, the
operator $\mathscr{A}$ in the space $L^\infty_{x_1}
L^2_{x_2,...,x_n} L^2_t (\mathbb{R}^{1+n})$ has succeed in absorbing
the partial derivative $\partial_{x_1}$. However,  it seem that
$\mathscr{A}$ can not deal with the partial derivative
$\partial_{x_2}$ in the space
 $L^\infty_{x_1} L^2_{x_2,...,x_n} L^2_t (\mathbb{R}^{1+n})$. So, we need a new way to handle the interaction
between $L^\infty_{x_1} L^2_{x_2,...,x_n} L^2_t (\mathbb{R}^{1+n})$
and $\partial_{x_2}$. We have the following

\begin{prop}\label{sm-eff-interact1}
Let $i=2,...,n$, $2\le q \le \infty$, $q>4/n$. Let $2\le r <\infty$,
$2/\gamma(r)= n(1/2-1/r)$, $\gamma\ge \gamma(r), \, \gamma>2.$ Then
we have
 \begin{align}
& \left \|\Box_k \partial_{x_i}  \mathscr{A}  f
\right\|_{L^\infty_{x_1} L^2_{x_2,...,x_n} L^2_t (\mathbb{R}^{1+n})}
\lesssim  \| \partial_{x_i} \partial^{-1}_{x_1} \Box_k
f\|_{L^1_{x_1} L^2_{x_2,...,x_n}L^2_t (\mathbb{R}^{1+n})},
\label{sm-int-1}\\
& \left \|\Box_k \partial_{x_i}  \mathscr{A}  f
\right\|_{L^\infty_{x_1} L^2_{x_2,...,x_n} L^2_t (\mathbb{R}^{1+n})}
\lesssim  \| \partial_{x_i} D^{-1/2}_{x_1} \Box_k f\|_{L^{\gamma'}
L^{r'}_x (\mathbb{R}^{1+n})},
\label{sm-int-1a}\\
& \left \|\Box_k \partial_{x_i}  \mathscr{A}  f \right\|_{L^q_{x_1}
L^\infty_{x_2,...,x_n} L^\infty_t (\mathbb{R}^{1+n})} \lesssim
\langle k_i\rangle^{1/2}  \langle k_1\rangle^{1/q} \|
 \Box_k  f\|_{L^1_{x_i}
L^2_{(x_j)_{j\not=i}}L^2_t (\mathbb{R}^{1+n})},
\label{sm-int-2}\\
& \left \|\Box_k \partial_{x_i}  \mathscr{A}  f \right\|_{L^q_{x_1}
L^\infty_{x_2,...,x_n} L^\infty_t (\mathbb{R}^{1+n})} \lesssim
\langle k_i\rangle   \langle k_1 \rangle^{1/q} \|
 \Box_k  f\|_{L^{\gamma'}_{t}
L^{r'}_{x}  (\mathbb{R}^{1+n})}. \label{sm-int-2a}
\end{align}
In \eqref{sm-int-2}, $q>2$ is required.
\end{prop}
{\bf Proof.} \eqref{sm-int-1} is a straightforward consequence of
Proposition \ref{GSE1}. We have
\begin{align}
 \mathcal {L} (\partial_{x_2} f, \psi) &:= \left|\int_\mathbb{R}\left(
\int_{\mathbb{R}}S(t-\tau) \partial_{x_2}
f(\tau) d\tau, \  \psi(t) \right ) dt \right| \nonumber\\
& \le \left\| \int_{\mathbb{R}}S(-\tau) \partial_{x_2}
D^{-1/2}_{x_1} f(\tau) d\tau\right\|_{L^2(\mathbb{R}^n)}
\left\|D^{1/2}_{x_1} \int_\mathbb{R} S(-t) \psi(t) dt
\right\|_{L^2(\mathbb{R}^n)}.
 \label{sm-int-3}
\end{align}
By the Strichartz inequality and Proposition \ref{GSE3},
\begin{align}
 \mathcal {L} (\partial_{x_2} f, \psi)
& \lesssim \| \partial_{x_2} D^{-1/2}_{x_1} f\|_{L^{\gamma'}_t
L^{r'}_x (\mathbb{R}^{1+n})}  \|\psi \|_{L^1_{x_1}
L^2_{x_2,...,x_n}L^2_t (\mathbb{R}^{1+n})} .
 \label{sm-int-3a}
\end{align}
By duality, \eqref{sm-int-3a} implies \eqref{sm-int-1a} in the case
$r>2$. In the case $r=2$, in view of the $1/2$-order smooth effect
of $S(t)$, we see that \eqref{sm-int-1a} also holds true. Similarly,
in view of Propositions \ref{GSE3}, \ref{MaxFunct-Mod-Dual} and
Lemma \ref{discret-deriv},
\begin{align}
 \mathcal {L} (\partial_{x_2} \Box_k f, \psi)
& \le \left\| \int_{\mathbb{R}}S(-\tau) D^1_{x_2} \Box_k f(\tau)
d\tau\right\|_{L^2(\mathbb{R}^n)} \left\|\tilde{\Box}_k
\int_\mathbb{R} S(-t) \psi(t) dt \right\|_{L^2(\mathbb{R}^n)}
\nonumber\\
& \lesssim \langle k_2\rangle^{1/2} \|  \Box_k f\|_{L^1_{x_2}
L^2_{(x_j)_{j\not=2}}L^2_t (\mathbb{R}^{1+n})} \langle k_1
\rangle^{1/q} \|\tilde{\Box}_k \psi \|_{L^{q'}_{x_1}
L^1_{x_2,...,x_n}L^1_t (\mathbb{R}^{1+n})} \nonumber\\
& \lesssim \langle k_2\rangle^{1/2}\langle k_1 \rangle^{1/q} \|
\Box_k f\|_{L^1_{x_2} L^2_{(x_j)_{j\not=2}}L^2_t (\mathbb{R}^{1+n})}
\|\psi \|_{L^{q'}_{x_1} L^1_{x_2,...,x_n}L^1_t (\mathbb{R}^{1+n})}.
 \label{sm-int-4}
\end{align}
By duality, \eqref{sm-int-2} follows from \eqref{sm-int-4}. Finally,
\begin{align}
 \mathcal {L} (\partial_{x_2} \Box_k f, \psi)
& \le \left\| \int_{\mathbb{R}}S(-\tau) \partial_{x_2} \Box_k
f(\tau) d\tau\right\|_{L^2(\mathbb{R}^n)} \left\|\tilde{\Box}_k
\int_\mathbb{R} S(-t) \psi(t) dt \right\|_{L^2(\mathbb{R}^n)}
\nonumber\\
& \lesssim \langle k_2\rangle  \langle k_1\rangle^{1/q}
\|\psi\|_{L^{q'}_{x_1} L^1_{x_2,...,x_n}L^1_t (\mathbb{R}^{1+n})} \|
 \Box_k
f \|_{L^{\gamma'}_{t} L^{r'}_{x} (\mathbb{R}^{1+n})}.
 \label{sm-int-4a}
\end{align}
If $r>2$ or $q>2$, \eqref{sm-int-4a} and Christ-Kiselev's Lemma
imply \eqref{sm-int-2a}, as desired. If $r=q=2$, in view of
Proposition \ref{MaxFunct-Mod}, we have also \eqref{sm-int-2a}.
$\hfill\Box$

\begin{lem}\label{sm-ef-int2}
Let $\psi : [0, \infty) \to [0,1]$ be a smooth bump function
satisfying $\psi(x)=1$ as $|x| \le 1$ and $\psi (x)=0$ if $|x| \ge
2$. Denote $\psi_1(\xi) = \psi(\xi_2/2\xi_1)$,  $\psi_2(\xi) =
1-\psi(\xi_2/2\xi_1)$, $\xi \in \mathbb{R}^n$. Then we have for
$\sigma \ge 0$,
 \begin{align}
& \sum_{k\in \mathbb{Z}^n, \, |k_1|>4} \langle k_1\rangle^\sigma
\left \|\mathscr{F}^{-1}_{\xi_1, \xi_2} \psi_1 \mathscr{F}_{x_1,
x_2} \Box_k
\partial_{x_2} \mathscr{A} f \right\|_{L^\infty_{x_1}
L^2_{x_2,...,x_n} L^2_t (\mathbb{R}^{1+n})}  \nonumber\\
& \quad \quad\quad  \lesssim \sum_{k\in \mathbb{Z}^n, \, |k_1|>4}
\langle k_1\rangle^\sigma \left \|\Box_k
 f \right\|_{L^1_{x_1}
L^2_{x_2,...,x_n} L^2_t (\mathbb{R}^{1+n})},  \label{sm-int-5a}
\end{align}
and for $\sigma \ge 1$,
 \begin{align}
& \sum_{k\in \mathbb{Z}^n, \, |k_1|>4} \langle k_1\rangle^\sigma
\left \|\mathscr{F}^{-1}_{\xi_1, \xi_2} \psi_2 \mathscr{F}_{x_1,
x_2} \Box_k
\partial_{x_2} \mathscr{A} f \right\|_{L^\infty_{x_1}
L^2_{x_2,...,x_n} L^2_t (\mathbb{R}^{1+n})}  \nonumber\\
 &  \quad \quad\quad \lesssim \sum_{k\in \mathbb{Z}^n, \,
|k_2|>4} \langle k_2\rangle^\sigma \left \|\Box_k
 f \right\|_{L^1_{x_1}
L^2_{x_2,...,x_n} L^2_t (\mathbb{R}^{1+n})}. \label{sm-int-5b}
\end{align}
\end{lem}
{\bf Proof.}  For simplicity, we denote
$$
I= \left \|\mathscr{F}^{-1}_{\xi_1, \xi_2} \psi_1 \mathscr{F}_{x_1,
x_2} \Box_k
\partial_{x_2} \mathscr{A} f \right\|_{L^\infty_{x_1}
L^2_{x_2,...,x_n} L^2_t (\mathbb{R}^{1+n})},
 $$
$$
II= \left \|\mathscr{F}^{-1}_{\xi_1, \xi_2} \psi_2 \mathscr{F}_{x_1,
x_2} \Box_k
\partial_{x_2} \mathscr{A} f \right\|_{L^\infty_{x_1}
L^2_{x_2,...,x_n} L^2_t (\mathbb{R}^{1+n})}.
$$

 Let $\eta_k $ be as in Lemma \ref{eq-mod-funct}. For $k \in
\mathbb{Z}^n$, $|k_1|
> 4$, applying the almost orthogonality of $\Box_k$,  we have
\begin{align}
I & \lesssim  \sum_{|\ell_1|, |\ell_2| \le 1}  \left
\|\mathscr{F}^{-1}_{\xi_1, \xi_2} \psi
\left(\frac{\xi_2}{2\xi_1}\right) \frac{\xi_2}{\xi_1}
\prod_{i=1,2}\eta_{k_i+\ell_i} (\xi_i) \mathscr{F}_{x_1, x_2} \Box_k
\partial_{x_1} \mathscr{A} f \right\|_{L^\infty_{x_1}
L^2_{x_2,...,x_n} L^2_t (\mathbb{R}^{1+n})}. \label{sm-int-7}
\end{align}
Denote
\begin{align}
(f \circledast_{12} g)(x) = \int_{\mathbb{R}^2}  f(t, x_1-y_1,
x_2-y_2, x_3,..., x_n) g(t, y_1,y_2) dy_1 dy_2. \label{sm-int-8}
\end{align}
We have for any Banach function space $X$ defined on
$\mathbb{R}^{1+n}$,
\begin{align}
\|f \circledast_{12} g\|_X \le \|g\|_{L^1_{y_1,y_2}(\mathbb{R}^2)}
\sup_{y_1, y_2} \| f(\cdot, \cdot-y_1, \cdot-y_2, \cdot,...,
\cdot)\|_X. \label{sm-int-9}
\end{align}
Hence, by \eqref{sm-int-7} and \eqref{sm-int-9},
\begin{align}
I & \lesssim  \sum_{|\ell_1|, |\ell_2| \le 1}  \left
\|\mathscr{F}^{-1}_{\xi_1, \xi_2} \psi
\left(\frac{\xi_2}{2\xi_1}\right) \frac{\xi_2}{\xi_1}
\prod_{i=1,2}\eta_{k_i+\ell_i} (\xi_i) \right\|_{L^1(\mathbb{R}^2)}
\left\|\Box_k
\partial_{x_1} \mathscr{A} f \right\|_{L^\infty_{x_1}
L^2_{x_2,...,x_n} L^2_t (\mathbb{R}^{1+n})}. \label{sm-int-10}
\end{align}
Using Bernstein's multiplier estimate, for $|k_1|>4$, we have
\begin{align}
& \!\! \!\! \left \|\mathscr{F}^{-1}_{\xi_1, \xi_2} \psi
\left(\frac{\xi_2}{2\xi_1}\right) \frac{\xi_2}{\xi_1}
\prod_{i=1,2}\eta_{k_i+\ell_i} (\xi_i) \right\|_{L^1(\mathbb{R}^2)}
\nonumber\\
& \lesssim  \sum_{|\alpha| \le 2} \left \|D^\alpha \left[\psi
\left(\frac{\xi_2}{2\xi_1}\right) \frac{\xi_2}{ \xi_1}
\prod_{i=1,2}\eta_{k_i+\ell_i} (\xi_i) \right]
\right\|_{L^2(\mathbb{R}^2)} \lesssim 1.
 \label{sm-int-11}
\end{align}
By Proposition \ref{1order-sm-mod}, \eqref{sm-int-10} and
\eqref{sm-int-11}, we have
\begin{align}
I & \lesssim \left\|\Box_k f \right\|_{L^1_{x_1} L^2_{x_2,...,x_n}
L^2_t (\mathbb{R}^{1+n})}, \ \ |k_1| \ge 4. \label{sm-int-12}
\end{align}
Next, we consider the estimate of $II$. Using Proposition
\ref{sm-eff-interact1},
\begin{align}
II & \lesssim \left \| \mathscr{F}^{-1}_{\xi_1, \xi_2}(\xi_2/\xi_1)
\psi_2 \mathscr{F}_{x_1,x_2}  \Box_k
 f \right\|_{L^1_{x_1}
L^2_{x_2,...,x_n} L^2_t (\mathbb{R}^{1+n})} \nonumber\\
& \lesssim \sum_{|\ell_1|, |\ell_2| \le 1}  \left
\|\mathscr{F}^{-1}_{\xi_1, \xi_2} \left(1-\psi \left(\frac{\xi_2}{2
\xi_1}\right)\right) \frac{\xi_2}{\xi_1}
\prod_{i=1,2}\eta_{k_i+\ell_i} (\xi_i) \right\|_{L^1(\mathbb{R}^2)}
\nonumber\\
& \quad \quad \quad \times  \left\|\Box_k
 f \right\|_{L^1_{x_1}
L^2_{x_2,...,x_n} L^2_t (\mathbb{R}^{1+n})}. \label{sm-int-13}
\end{align}
Notice that ${\rm supp }\psi_2 \subset \{\xi: \, |\xi_2| \ge
2|\xi_1|\}$. If $|k_1| \ge 4$, we have $|k_2|>6$ and $|k_2| \ge
|k_1|$ in the summation of the left-hand side of \eqref{sm-int-5b}.
So, $\sum_{k\in \mathbb{Z}^n, \ |k_1|>4}\langle k_1\rangle^{\sigma}
II \le \sum_{k\in \mathbb{Z}^n, \ |k_1|>4}\langle
k_2\rangle^{\sigma-1} \langle k_1\rangle  II$.
\begin{align}
& \left \|\mathscr{F}^{-1}_{\xi_1, \xi_2} \left(1-\psi
\left(\frac{\xi_2}{2 \xi_1}\right)\right) \frac{\xi_2}{\xi_1}
\prod_{i=1,2}\eta_{k_i+\ell_i} (\xi_i) \right\|_{L^1(\mathbb{R}^2)} \nonumber\\
& \lesssim  \sum_{|\alpha|\le 2} \left \|D^\alpha
\left[\mathscr{F}^{-1}_{\xi_1, \xi_2} \left(1-\psi
\left(\frac{\xi_2}{2 \xi_1}\right)\right) \frac{\xi_2}{ \xi_1}
\prod_{i=1,2}\eta_{k_i+\ell_i} (\xi_i)\right]
\right\|_{L^2(\mathbb{R}^2)} \nonumber\\
&  \lesssim \langle k_2 \rangle  \langle k_1\rangle^{-1}.
 \label{sm-int-14}
\end{align}
\eqref{sm-int-13} and \eqref{sm-int-14} yield the estimate of $II$,
as desired. $\hfill \Box$

\begin{con}\label{rem-int}
\rm Using a similar way as in the proof of \eqref{sm-int-1}, we can
show that
\begin{align}
 \left \|\Box_k \partial_{x_2}  \int_{\mathbb{R}} S(t-\tau) f(\tau)
d\tau \right\|_{L^\infty_{x_1} L^2_{x_2,...,x_n} L^2_t
(\mathbb{R}^{1+n})}  \lesssim  \| D^{1/2}_{x_2} D^{-1/2}_{x_1}
\Box_k f\|_{L^1_{x_2} L^2_{x_1,x_3,...,x_n}L^2_t
(\mathbb{R}^{1+n})}. \nonumber
\end{align}
So, we can conjecture that
\begin{align}
& \left \|\Box_k \partial_{x_2} \mathscr{A} f
\right\|_{L^\infty_{x_1} L^2_{x_2,...,x_n} L^2_t (\mathbb{R}^{1+n})}
 \lesssim  \| D^{1/2}_{x_2} D^{-1/2}_{x_1} \Box_k f\|_{L^1_{x_2}
L^2_{x_1,x_3,...,x_n}L^2_t (\mathbb{R}^{1+n})}. \label{sm-int-rem1}
\end{align}
Since we do not know if the Christ-Kiselev Lemma holds in the
endpoint case, it is not clear for us if \eqref{sm-int-rem1} is
true.

If \eqref{sm-int-rem1} is true, repeating the proof above, we can
show that \eqref{sm-int-5a} holds for all $\sigma \ge 1/2$. We can
improve the results of Theorems \ref{DNLS-mod} and \ref{DNLS1-mod}
by assuming that $u_0 \in M^{1+1/m}_{2,1}$ and $u_0 \in
M^{1/m}_{2,1}$, respectively.

\end{con}

\begin{lem}\label{sm-ef-int3}
Let  $2\le q \le \infty$, $q>4/n$ and $(\gamma, \ r)$ be as in
Proposition \ref{sm-eff-interact1}. Let $k=(k_1,...,k_n)$, $k_{\rm
max}:=\max_{1\le i \le n} |k_i|$.  Then we have
 \begin{align}
& \left \|\Box_k \partial_{x_i}  \mathscr{A}  f \right\|_{L^q_{x_1}
L^\infty_{x_2,...,x_n} L^\infty_t (\mathbb{R}^{1+n})} \lesssim
\langle k_{\max}\rangle^{1+1/q} \|
 \Box_k  f\|_{L^{\gamma'}_{t}
L^{r'}_{x}  (\mathbb{R}^{1+n})}. \label{sm-int-i1a}
\end{align}
\end{lem}
{\bf Proof.}   It follows from \eqref{sm-int-2a} that
\eqref{sm-int-i1a} holds.   $\hfill \Box$

\begin{lem}\label{sm-ef-int4}
Let $k=(k_1,...,k_n)$, $k_{\rm max}:=\max_{1\le i \le n} |k_i|$ and
$q> 2\vee 4/n.$ Then we have for $\sigma \ge 0$ and $i,
\alpha=1,...,n$,
 \begin{align}
& \sum_{k\in \mathbb{Z}^n, \, |k_\alpha|=k_{ \max}>4} \langle k
\rangle^\sigma \left \|\Box_k
\partial_{x_i} \mathscr{A} f \right\|_{L^q_{x_1}
L^\infty_{x_2,...,x_n} L^\infty_t (\mathbb{R}^{1+n})}  \nonumber\\
 &  \quad \quad\quad \lesssim \sum_{k\in \mathbb{Z}^n, \,
|k_\alpha|>4} \langle k_\alpha \rangle^{\sigma+ 1/2+1/q} \left
\|\Box_k
 f \right\|_{L^1_{x_\alpha}
L^2_{(x_j)_{j\not=\alpha}} L^2_t (\mathbb{R}^{1+n})}.
\label{sm-int-4.4a}
\end{align}
\end{lem}
{\bf Proof.} First, we consider the case $\alpha=1$. In view of
\eqref{st-sm-m-c6} and $|k_1|= k_{\max}>4$,
\begin{align}
  \left \|\Box_k
\partial_{x_i} \mathscr{A} f \right\|_{L^{q}_{x_1}
L^\infty_{x_2,...,x_n} L^\infty_t (\mathbb{R}^{1+n})}  & \ \lesssim
\sum_{|\ell_1|, |\ell_i| \le 1} \left\| \mathscr{F}^{-1}_{\xi_1,
\xi_i}  \left(\frac{\xi_i}{\xi_1} \eta_{k_i+\ell_i} (\xi_i)
\eta_{k_1+\ell_1} (\xi_1) \right)\right\|_{L^1(\mathbb{R}^2)}
\nonumber\\
 & \ \ \ \   \times \| \Box_k
\partial_{x_1} \mathscr{A} f \|_{L^{q}_{x_1}
L^\infty_{x_2,...,x_n} L^\infty_t (\mathbb{R}^{1+n})} \nonumber\\
 & \ \lesssim \langle k_i\rangle \langle k_1\rangle^{-1} \langle k_1\rangle^{1/2+1/q} \| \Box_k
 f \|_{L^{1}_{x_1} L^2_{x_2,...,x_n} L^2_t
(\mathbb{R}^{1+n})} \nonumber\\
 & \ \lesssim  \langle k_1\rangle^{1/2+1/q} \| \Box_k
 f \|_{L^{1}_{x_1} L^2_{x_2,...,x_n} L^2_t
(\mathbb{R}^{1+n})} . \label{sm-int-4.4b}
\end{align}
\eqref{sm-int-4.4b} implies the result, as desired. Next, we
consider the case $\alpha=2$. Notice that $|k_2| =\max_{1\le i \le
n} |k_i|>4$. By \eqref{sm-int-2},
\begin{align}
  \left \|\Box_k
\partial_{x_i} \mathscr{A} f \right\|_{L^{q}_{x_1}
L^\infty_{x_2,...,x_n} L^\infty_t (\mathbb{R}^{1+n})}  & \ \lesssim
\sum_{|\ell_2|, |\ell_i| \le 1} \left\| \mathscr{F}^{-1}_{\xi_i,
\xi_2}  \left(\frac{\xi_i}{\xi_2} \eta_{k_i+\ell_i} (\xi_i)
\eta_{k_2+\ell_2} (\xi_2) \right)\right\|_{L^1(\mathbb{R}^2)}
\nonumber\\
 & \ \ \ \   \times \| \Box_k
\partial_{x_2} \mathscr{A} f \|_{L^{q}_{x_1}
L^\infty_{x_2,...,x_n} L^\infty_t (\mathbb{R}^{1+n})} \nonumber\\
 & \ \lesssim \langle k_i\rangle \langle k_2\rangle^{-1} \langle k_2\rangle^{1/2+1/q} \| \Box_k
 f \|_{L^{1}_{x_2} L^2_{x_1,x_3,...,x_n} L^2_t
(\mathbb{R}^{1+n})} \nonumber\\
 & \ \lesssim  \langle k_2\rangle^{1/2+1/q} \| \Box_k
 f \|_{L^{1}_{x_2} L^2_{x_1,x_3,...,x_n} L^2_t
(\mathbb{R}^{1+n})} . \label{sm-int-4.4c}
\end{align}
The other cases $\alpha=3,...,n$ is analogous to the case $\alpha=2$
and we omit the details of the proof. $\hfill \Box$

\begin{rem} \rm
From the proof of Lemma \ref{sm-ef-int4}, we easily see that
\begin{align}
& \sum_{k\in \mathbb{Z}^n, \, |k_\alpha|=k_{ \max}>4} \langle k
\rangle^\sigma \left \|\Box_k
\partial_{x_i} \mathscr{A} f \right\|_{L^q_{x_\beta}
L^\infty_{(x_j)_{j\not=\beta}} L^\infty_t (\mathbb{R}^{1+n})}  \nonumber\\
 &  \quad \quad\quad \lesssim \sum_{k\in \mathbb{Z}^n, \,
|k_\alpha|>4} \langle k_\alpha \rangle^{\sigma+ 1/2+1/q} \left
\|\Box_k
 f \right\|_{L^1_{x_\alpha}
L^2_{(x_j)_{j\not=\alpha}} L^2_t (\mathbb{R}^{1+n})}.
\label{sm-int-4.4d}
\end{align}
\end{rem}

\section{Proof of Theorem \ref{DNLS1-mod}} \label{pf-thm1}
Now we briefly indicate the proof of Theorem \ref{DNLS1-mod}.  We
assume that the nonlinear term takes the form
\begin{align*}
F(u, \nabla u) = \partial_{x_1} (u^{\kappa_1+1})+ \partial_{x_2}
(u^{\kappa_2+1}).
\end{align*}
In order to handle the nonlinear term $\partial_{x_i}
(u^{\kappa_i+1})$, we use the space $L^\infty_{x_i}
L^2_{(x_j)_{j\not=i}}L^2_t (\mathbb{R}^{1+n})$ to absorb the
derivative $\partial_{x_i}$. Hence, we introduce the following
semi-norms to treat the nonlinearity:
\begin{align*}
& \| u\|_{Y_i} = \sum_{k\in \mathbb{Z}^n, \ |k_i|>4} \langle
k_i\rangle  \|\Box_k u\|_{L^{\infty}_{x_i} L^2_{(x_j)_{j\not=
i}}L^2_t(\mathbb{R}^{1+n})}, \quad i=1,2.
\end{align*}
Since \eqref{1-sm-mod} is a worse estimate in the case $|k_i|
\lesssim 1$,  we throw away the low frequency part in the
$\xi_i$-direction in the definition of $\|u\|_{Y_i}$. To handle the
low frequency part, we use the Strichartz norm:
\begin{align*}
\|u\|_S = \sum_{k\in \mathbb{Z}^n} \langle k\rangle^{1/2} \|\Box_k
u\|_{L^\infty_tL^2_x \bigcap L^{2+\kappa}_{x, t}
(\mathbb{R}^{1+n})}.
\end{align*}
We emphasize that the Strichartz inequalities \eqref{st-sm-m-c1} and
\eqref{st-sm-m-c9} are better estimates than the smooth effects in
\eqref{st-sm-m-c2} and \eqref{st-sm-mo-5} for the low frequency
part, respectively. Using the integral equation
\begin{align*}
 u(t) = S(t) u_0 - {\rm i} \mathscr{A} (\partial_{x_1}
u^{\kappa_1+1} + \partial_{x_2} u^{\kappa_2+1} ),
\end{align*}
we have
\begin{align*}
\| u\|_{Y_1} \le  \|S(t) u_0\|_{Y_1} +  \| \mathscr{A}
(\partial_{x_1} u^{\kappa_1+1})\|_{Y_1} + \|\mathscr{A}
(\partial_{x_2} u^{\kappa_2+1} )\|_{Y_1}.
\end{align*}
In view of \eqref{st-sm-m-c2},  $\|S(t) u_0\|_{Y_1}$ is bounded by
$\|u_0\|_{M^{1/2}_{2,1}}$. $\| \mathscr{A} (\partial_{x_1}
u^{\kappa_1+1})\|_{Y_1} $ can be handled by using the linear
estimates obtained in Section \ref{LE1}.  Noticing that
$$
u^{\kappa_1+1}= \left(\sum_{k^{(1)},..., k^{(\kappa_1+1)} \in
\mathbb{S}^1} + \sum_{k^{(1)},..., k^{(\kappa_1+1)} \in \mathbb{Z}^n
\setminus \mathbb{S}^1} \right) \Box_{k^{(1)}}u ...
\Box_{k^{(\kappa_1+1)}}u,
$$
 where $\mathbb{S}^1= \{k^{(1)},...,
k^{(\kappa_1+1)} \in \mathbb{Z}^n: \ |k^{(1)}_1|\vee ...\vee
|k^{(\kappa_1+1)}_1| >4\}$, \eqref{st-sm-mo-5} and
\eqref{st-sm-mo-7} in Corollary \ref{st-sm-m-c} yield,
\begin{align}
\|\mathscr{A} (\partial_{x_1} u^{\kappa_1+1})\|_{Y_1} &  \lesssim
\!\!\!\! \sum_{k\in \mathbb{Z}^n, \ |k_1|>4} \langle k_1\rangle
\!\! \sum_{\mathbb{S}^1} \|\Box_k (\Box_{k^{(1)}}u ...
\Box_{k^{(\kappa_1+1)}}u )\|_{L^1_{x_1}
L^2_{x_2,..., x_n}L^2_t} \nonumber\\
&  \ \ \  + \sum_{k\in \mathbb{Z}^n, \ |k_1|> 4} \langle
k_1\rangle^{3/2} \!\! \sum_{\mathbb{Z}^n \setminus \mathbb{S}^1}
\|\Box_k (\Box_{k^{(1)}}u ... \Box_{k^{(\kappa_1+1)}}u
)\|_{L^{\frac{2+\kappa}{1+\kappa}}_{x, t}} . \label{control2}
\end{align}
By performing a nonlinear mapping estimate, we have
\begin{align}
\|\mathscr{A} (\partial_{x_1} u^{\kappa_1+1})\|_{Y_1} &  \lesssim
 \|u\|_{Y_1} \|u\|^{\kappa_1}_{Z_1}+ \|u\|^{\kappa_1+1}_S,
 \label{control3}
\end{align}
where
$$
\|u\|_{Z_i} = \sum_{k \in \mathbb{Z}^n} \langle
k\rangle^{1/2-1/\kappa} \|\Box_k u\|_{L^\kappa_{x_i}
L^\infty_{(x_j)_{j\not=i}} L^\infty_t(\mathbb{R}^{1+n}) }, \ \
i=1,2.
$$

Unfortunately, $\|\mathscr{A} (\partial_{x_2} u^{\kappa_2+1}
)\|_{Y_1}$ contains the interaction between the working space
$L^\infty_{x_1} L^2_{x_2,...,x_n}L^2_t (\mathbb{R}^{1+n})$ and the
derivative $\partial_{x_2}$, which is out of the control of the
smooth effect \eqref{st-sm-mo-5}. So, we look for another way to
estimate $\|\mathscr{A} (\partial_{x_2} u^{\kappa_2+1} )\|_{Y_1}$.
Roughly speaking, our idea is to use the following estimates (see
Lemma \ref{sm-ef-int2}):
\begin{align}
\|\mathscr{F}^{-1}\chi_{\{\xi: \ |\xi_2|\le |\xi_1|\}}\mathscr{F}
\mathscr{A} (\partial_{x_2} f )\|_{Y_1} & \lesssim \sum_{k\in
\mathbb{Z}^n \; |k_1|>4}  \langle k_1\rangle  \|\Box_k
f\|_{L^1_{x_1} L^2_{x_2,...,x_n}L^2_t (\mathbb{R}^{1+n})}, \nonumber
\\
\|\mathscr{F}^{-1}\chi_{\{\xi: \ |\xi_1|\le |\xi_2|\}}\mathscr{F}
\mathscr{A} (\partial_{x_2} f )\|_{Y_1} & \lesssim \sum_{k\in
\mathbb{Z}^n \; |k_2|>4}  \langle k_2\rangle  \|\Box_k
f\|_{L^1_{x_1} L^2_{x_2,...,x_n}L^2_t (\mathbb{R}^{1+n})}, \nonumber
\end{align}
where $\chi_E$ denotes the characteristic function on the set $E$.
So, $\|\mathscr{A} (\partial_{x_2} u^{\kappa_2+1} )\|_{Y_1}$ has
similar bound to $\|\mathscr{A} (\partial_{x_1} u^{\kappa_2+1}
)\|_{Y_1}$ as in \eqref{control2}. Eventually, we have
\begin{align}
\|\mathscr{A} (\partial_{x_2} u^{\kappa_2+1})\|_{Y_1} &  \lesssim
 (\|u\|_{Y_1}+\|u\|_{Y_2}) \|u\|^{\kappa_2}_{Z_1}+
 \|u\|^{\kappa_2+1}_S. \label{control5}
\end{align}
By using the integral equation, we need to further bound $
\|\mathscr{A}
\partial_{x_i} u^{\kappa_i+1}\|_{Z_1 \bigcap S} $, $i=1,2$. For
instance,
 for the estimate of $
\|\mathscr{A}
\partial_{x_2} u^{\kappa_2+1}\|_{Z_1} $, we resort to the above idea
and consider the following interaction estimate:
\begin{align*}
& \left \|\Box_k \partial_{x_2}  \mathscr{A}  f
\right\|_{L^\kappa_{x_1} L^\infty_{x_2,...,x_n} L^\infty_t
(\mathbb{R}^{1+n})} \lesssim \langle k_2\rangle^{1/2}\langle k
\rangle^{1/\kappa} \|\Box_k f\|_{L^1_{x_2}
L^2_{(x_j)_{j\not=2}}L^2_t (\mathbb{R}^{1+n})},
\end{align*}
which leads to that we can bound $ \|\mathscr{A} \partial_{x_2}
u^{\kappa_2+1}\|_{Z_1}$  by an analogous version of the right-hand
side of \eqref{control2}, so, by \eqref{control3} (see Lemma
\ref{sm-ef-int4}).

Finally, using \eqref{st-sm-m-c4} and \eqref{st-sm-m-c9},  we can
get the same estimate of $\|\mathscr{A}
\partial_{x_i} u^{\kappa_i+1}\|_{S} $ as in \eqref{control3} and
\eqref{control5}, respectively.

\medskip

{\bf Proof of Theorem \ref{DNLS1-mod}.} We now give the details of
the proof of Theorem \ref{DNLS1-mod}.  Denote
\begin{align*}
& \rho_1(u) = \sum^n_{i=1}\sum_{k\in \mathbb{Z}^n, \ |k_i|>4}
\langle k_i\rangle  \|\Box_k u\|_{L^{\infty}_{x_i}
L^2_{(x_j)_{j\not=i}}L^2_t(\mathbb{R}^{1+n})},\\
& \rho_2(u) = \sum^n_{i=1}\sum_{k\in \mathbb{Z}^n} \langle
k\rangle^{1/2-1/\kappa} \|\Box_k u\|_{L^{\kappa}_{x_i}
L^\infty_{(x_j)_{j\not=i}}L^\infty_t(\mathbb{R}^{1+n})},\\
& \rho_3(u) = \sum_{k\in \mathbb{Z}^n} \langle k\rangle^{1/2}
\|\Box_k u\|_{L^\infty_tL^2_x  \bigcap L^{2+\kappa}_{x, t}
(\mathbb{R}^{1+n})}.
\end{align*}
Put
\begin{align*}
X:= \left\{u\in \mathscr{S}'(\mathbb{R}^{1+n}): \ \|u\|_X:=
\sum^3_{i=1}\rho_i(u) \le \delta_0 \right\}.
\end{align*}
We consider the following mapping:
\begin{align*}
\mathscr{T}: u(t) \to S(t) u_0 - {\rm i} \mathscr{A}
\left(\sum^n_{i=1} \lambda_i \partial_{x_i} u^{\kappa_i+1} \right).
\end{align*}
For convenience, we denote
\begin{align*}
& \|u\|_{Y_i} = \sum_{k\in \mathbb{Z}^n, \ |k_i|>4} \langle
k_i\rangle  \|\Box_k u\|_{L^{\infty}_{x_i}
L^2_{(x_j)_{j\not=i}}L^2_t(\mathbb{R}^{1+n})}.
\end{align*}
In order to estimate $\rho_1(u)$, it suffices to control
$\|\cdot\|_{Y_1}$.  By \eqref{smo-eff.21} and Plancherel's identity,
we have
\begin{align*}
\|S(t) u_0\|_{Y_1} & \lesssim  \sum_{k\in \mathbb{Z}^n, \ |k_1|>4}
\langle k_1\rangle  \|\Box_k D^{-1/2}_{x_1}
u_0\|_{L^2(\mathbb{R}^{n})} \\
&  \lesssim \sum_{k\in \mathbb{Z}^n} \langle k_1\rangle^{1/2}
\|\Box_k u_0\|_{L^2(\mathbb{R}^{n})}.
\end{align*}
By \eqref{MFM-2}, Lemma \ref{Strichartz-mod}, we have
\begin{align*}
\rho_i (S(t) u_0) &  \lesssim \sum_{k\in \mathbb{Z}^n} \langle k
\rangle^{1/2} \|\Box_k u_0\|_{L^2(\mathbb{R}^{n})}, \ \ i=2,3.
\end{align*}
Denote
\begin{align*}
& \mathbb{S}^{(i)}_{\ell, 1}:= \{(k^{(1)},..., k^{(\kappa_\ell+1)})
\in (\mathbb{Z}^n)^n: \
|k^{(1)}_i|\vee...\vee |k^{(\kappa_\ell+1)}_i|>4\}, \\
& \mathbb{S}^{(i)}_{\ell, 2}:= \{(k^{(1)},..., k^{(\kappa_\ell+1)})
\in (\mathbb{Z}^n)^n: \ |k^{(1)}_i|\vee...\vee
|k^{(\kappa_\ell+1)}_i|\le 4\}.
\end{align*}
Using the frequency-uniform decomposition, we have
\begin{align}
u^{\kappa_\ell+1}  =  & \sum_{k^{(1)},..., k^{(\kappa_\ell+1)}\in
\mathbb{Z}^n} \Box_{k^{(1)}} u ... \Box_{k^{(\kappa_\ell+1)}} u  \nonumber\\
 = & \sum_{\mathbb{S}^{(i)}_{\ell,1}} \Box_{k^{(1)}} u ... \Box_{k^{(\kappa_\ell+1)}} u  +
\sum_{\mathbb{S}^{(i)}_{\ell,2}} \Box_{k^{(1)}} u ...
\Box_{k^{(\kappa_\ell+1)}} u . \label{decomp}
\end{align}
Using \eqref{st-sm-mo-5} and \eqref{st-sm-mo-7}, we obtain that
\begin{align}
\|\mathscr{A} \partial_{x_1} u^{\kappa_1+1}\|_{Y_1} & \lesssim
\sum_{k \in \mathbb{Z}^n, \ |k_1| >4} \langle k_1\rangle
\sum_{\mathbb{S}^{(1)}_{1,1}} \|\Box_k \left(\Box_{k^{(1)}} u ...
\Box_{k^{(\kappa_1+1)}} u \right) \|_{L^{1}_{x_1}
L^2_{x_2,...,x_n}L^2_t(\mathbb{R}^{1+n})}  \nonumber\\
& \ \  +  \sum_{k \in \mathbb{Z}^n, \ |k_1| >4} \langle
k_1\rangle^{3/2} \sum_{\mathbb{S}^{(1)}_{1,2}} \| \Box_k
\left(\Box_{k^{(1)}} u ... \Box_{k^{(\kappa_1+1)}} u
\right)\|_{L^{(2+\kappa)/(1+\kappa)}_{t,x}(\mathbb{R}^{1+n}) } \nonumber\\
& := I +II . \label{nonhomo-1}
\end{align}
In view of the support property of $\widehat{\Box_k u}$, we see that
\begin{align}
 \Box_k \left(\Box_{k^{(1)}} u ... \Box_{k^{(\kappa_1+1)}} u
 \right)=0, \ \ if  \ \ |k- k^{(1)}-...-k^{(\kappa_1+1)}| \ge C.
 \label{orth}
\end{align}
Hence, by Lemma \ref{discret-deriv},
\begin{align}
I & \lesssim \sum_{k \in \mathbb{Z}^n, \ |k_1| >4} \langle
k_1\rangle   \sum_{\mathbb{S}^{(1)}_{1,1}} \|\Box_{k^{(1)}} u ...
\Box_{k^{(\kappa_1+1)}} u  \|_{L^{1}_{x_1}
L^2_{x_2,...,x_n}L^2_t(\mathbb{R}^{1+n})} \chi_{|k-
k^{(1)}-...-k^{(\kappa_1+1)}| \le C}. \label{non-est-1}
\end{align}
By H\"older's inequality and $\|\Box_k u \|_{L^{\infty}_{x}}
\lesssim \|\Box_k u \|_{L^{2}_{x}} $ uniformly holds for all $k\in
\mathbb{Z}^n$, we have
\begin{align*}
& \|\Box_{k^{(1)}} u ... \Box_{k^{(\kappa_1+1)}} u \|_{L^{1}_{x_1}
L^2_{x_2,...,x_n}L^2_t(\mathbb{R}^{1+n})} \\
& \le \|\Box_{k^{(1)}} u\|_{L^{\infty}_{x_1}
L^2_{x_2,...,x_n}L^2_t(\mathbb{R}^{1+n})} \prod^{\kappa_1+1}_{i=2}
\|\Box_{k^{(i)}} u \|_{L^{\kappa}_{x_1}
L^\infty_{x_2,...,x_n}L^\infty_t \, \bigcap \, L^\infty_t
L^2_x(\mathbb{R}^{1+n})}.
\end{align*}
Since $|k- k^{(1)}-...-k^{(\kappa_1+1)}| \le C$ implies that $|k_1-
k^{(1)}_1-...-k^{(\kappa_1+1)}_1| \le C$, we see that $|k_1| \le C
\max_{i=1,...,\kappa_1+1} |k^{(i)}_1|$. We may assume that
$|k^{(1)}_1|= \max_{i=1,...,\kappa_1+1} |k^{(i)}_1|$ in the
summation $\sum_{\mathbb{S}^{(1)}_{1,1}}$ in \eqref{non-est-1}
above. So,
\begin{align}
I & \lesssim    \sum_{k^{(1)} \in \mathbb{Z}^n, \ |k^{(1)}_1|>4}
\langle k^{(1)}_1\rangle  \|\Box_{k^{(1)}}
u\|_{L^{\infty}_{x_1} L^2_{x_2,...,x_n}L^2_t(\mathbb{R}^{1+n})} \nonumber\\
& \ \ \times  \sum_{k^{(2)},...,k^{(\kappa_1+1)}\in \mathbb{Z}^n}
\prod^{\kappa_1+1}_{i=2} \|\Box_{k^{(i)}} u \|_{L^{\kappa}_{x_1}
L^\infty_{x_2,...,x_n}L^\infty_t \, \bigcap \, L^\infty_t
L^2_x (\mathbb{R}^{1+n})} \nonumber\\
& \lesssim \rho_1(u) (\rho_2(u)+\rho_3(u))^{\kappa_1} .
\label{est-I}
\end{align}
In view of \eqref{orth} we easily see that $|k_1| \le C$ in $II$ of
\eqref{nonhomo-1}. Hence,
\begin{align}
II & \lesssim \sum_{k \in \mathbb{Z}^n, \ |k_1| >4}
\sum_{\mathbb{S}^{(1)}_{1,2}} \|\Box_{k^{(1)}} u ...
\Box_{k^{(\kappa_1+1)}} u  \|_{L^{(2+\kappa)/(1+\kappa)}_{x, t}
(\mathbb{R}^{1+n})}
\chi_{|k- k^{(1)}-...-k^{(\kappa_1+1)}| \le C} \nonumber\\
& \lesssim  \sum_{\mathbb{S}^{(1)}_{1,2}} \|\Box_{k^{(1)}} u ...
\Box_{k^{(\kappa_1+1)}} u \|_{L^{(2+\kappa)/(1+\kappa)}_{x, t}
(\mathbb{R}^{1+n})} \nonumber\\
& \lesssim  \sum_{\mathbb{S}^{(1)}_{1,2}}
\prod^{\kappa_1+1}_{i=1}\|\Box_{k^{(i)}} u \|_{L^{2+\kappa}_{x, t}
\, \bigcap \, L^\infty_t L^2_x (\mathbb{R}^{1+n})} \lesssim
\rho_3(u)^{1+\kappa_1}. \label{non-est-2}
\end{align}
Hence, we have
\begin{align}
\|\mathscr{A} \partial_{x_1} u^{\kappa_1+1}\|_{Y_1} & \lesssim
\rho_1(u) (\rho_2(u)+\rho_3(u))^{\kappa_1} + \rho_3(u)^{1+\kappa_1}.
\label{est-lam1}
\end{align}

Next, we estimate $\|\mathscr{A} \partial_{x_2}
u^{\kappa_2+1}\|_{Y_1}$. Let $\psi_i$ be as in Lemma
\ref{sm-ef-int2}.  For convenience, we write
\begin{align}
P_i= \mathscr{F}^{-1}_{\xi_1, \xi_2} \psi_i \mathscr{F}_{x_1, x_2},
\ \ i=1,2. \label{Proj-i}
\end{align}
We have
 \begin{align}
\|\mathscr{A} \partial_{x_2} u^{\kappa_2+1}\|_{Y_1} \lesssim & \left
\|P_1
\partial_{x_2} \mathscr{A} u^{\kappa_2+1} \right\|_{Y_1} + \left \|P_2
\partial_{x_2} \mathscr{A} u^{\kappa_2+1} \right\|_{Y_1}:=III+IV. \label{est-int-1}
\end{align}
Using the decomposition \eqref{decomp},
\begin{align}
III & \le \Big \|P_1
\partial_{x_2} \mathscr{A} \sum_{\mathbb{S}^{(1)}_{2,1}} (\Box_{k^{(1)}} u ...
\Box_{k^{(\kappa_2+1)}} u )  \Big\|_{Y_1} \nonumber\\
& \quad  +  \Big \|P_1
\partial_{x_2} \mathscr{A} \sum_{\mathbb{S}^{(1)}_{2,2}} (\Box_{k^{(1)}} u ...
\Box_{k^{(\kappa_2+1)}} u )  \Big\|_{Y_1}:=III_1+III_2.
\label{est-int-1a}
\end{align}
Applying  Lemma \ref{sm-ef-int2} and then following the same way as
in the estimate to \eqref{non-est-1},
\begin{align}
III_1 & \lesssim \sum_{k \in \mathbb{Z}^n, \ |k_1| >4} \langle k_1
\rangle  \sum_{\mathbb{S}^{(1)}_{2,1}} \|\Box_k \left(\Box_{k^{(1)}}
u ... \Box_{k^{(\kappa_2+1)}} u \right) \|_{L^{1}_{x_1}
L^2_{x_2,...,x_n}L^2_t(\mathbb{R}^{1+n})}
\nonumber\\
&\lesssim \rho_1(u) (\rho_2(u)+\rho_3(u))^{\kappa_2}.
\label{est-int-1b}
\end{align}
For the estimate of $III_2$, noticing the fact that ${\rm supp}
\psi_1 \subset \{\xi: \ |\xi_2| \le 4 |\xi_1|\}$ and using the
multiplier estimate,  then applying \eqref{sm-int-1a}, we have
\begin{align}
III_2 & \lesssim   \sum_{k \in \mathbb{Z}^n, \ |k_1| >4, \ |k_2|
\lesssim |k_1|} \langle k_1\rangle^{3/2}
\sum_{\mathbb{S}^{(1)}_{2,2}} \|\Box_k \left(\Box_{k^{(1)}} u ...
\Box_{k^{(\kappa_2+1)}} u
\right)\|_{L^{(2+\kappa)/(1+\kappa)}_{t,x}(\mathbb{R}^{1+n})
}\nonumber\\
& \lesssim  \rho_3(u)^{1+\kappa_2} . \label{est-int-2}
\end{align}
We need to further control $IV$. Using the decomposition
\eqref{decomp},
\begin{align}
IV & \le \Big \|P_2
\partial_{x_2} \mathscr{A} \sum_{\mathbb{S}^{(2)}_{2,1}} (\Box_{k^{(1)}} u ...
\Box_{k^{(\kappa_2+1)}} u )  \Big\|_{Y_1} \nonumber\\
& \quad  +  \Big \|P_2
\partial_{x_2} \mathscr{A} \sum_{\mathbb{S}^{(2)}_{2,2}} (\Box_{k^{(1)}} u ...
\Box_{k^{(\kappa_2+1)}} u )  \Big\|_{Y_1}:=IV_1+IV_2.
\label{est-int-2a}
\end{align}
By Lemma \ref{sm-ef-int2},
\begin{align}
IV_1 & \lesssim \sum_{k \in \mathbb{Z}^n, \ |k_2| >4} \langle k_2
\rangle  \sum_{\mathbb{S}^{(2)}_{2,1}} \|\Box_k \left(\Box_{k^{(1)}}
u ... \Box_{k^{(\kappa_2+1)}} u \right) \|_{L^{1}_{x_1}
L^2_{x_2,...,x_n}L^2_t(\mathbb{R}^{1+n})}. \label{est-int-2b}
\end{align}
By symmetry of $k^{(1)},..., k^{(\kappa_2+1)}$, we can assume that
$|k^{(1)}_2 |=\max_{1\le i \le \kappa_2+1} |k^{(i)}_2 |$ in
$\mathbb{S}^{(2)}_{2,1}$. Using the same way as in the estimate of
$I$, we have
\begin{align}
IV_1 & \lesssim \sum_{\mathbb{S}^{(2)}_{2,1}, \ |k^{(1)}_2| >4}
\langle k^{(1)}_2 \rangle  \| \Box_{k^{(1)}} u ...
\Box_{k^{(\kappa_2+1)}} u
 \|_{L^{1}_{x_1} L^2_{x_2,...,x_n}L^2_t(\mathbb{R}^{1+n})}.
\label{est-int-2ba}
\end{align}
By H\"older's inequality,
\begin{align}
&  \| \Box_{k^{(1)}} u ... \Box_{k^{(\kappa_2+1)}} u
 \|_{L^{1}_{x_1} L^2_{x_2,...,x_n}L^2_t(\mathbb{R}^{1+n})}
 \nonumber\\
& \lesssim  \| \Box_{k^{(1)}} u |\Box_{k^{(2)}} u ...
\Box_{k^{(\kappa_2+1)}} u|^{1/2}
 \|_{L^2_{x,t} (\mathbb{R}^{1+n})} \nonumber\\
  &  \ \ \ \ \times  \||\Box_{k^{(2)}} u ...
\Box_{k^{(\kappa_2+1)}} u|^{1/2}
 \|_{L^{2}_{x_1} L^\infty_{x_2,...,x_n}L^\infty_t(\mathbb{R}^{1+n})} \nonumber\\
& \lesssim  \| \Box_{k^{(1)}} u \|_{L^{\infty}_{x_2}
L^2_{x_1,x_3,...,x_n}L^2_t(\mathbb{R}^{1+n})}
\prod^{\kappa_2+1}_{i=2} \|\Box_{k^{(i)}} u
 \|_{L^{\kappa_2}_{x_2} L^\infty_{x_1,x_3,...,x_n}L^\infty_t(\mathbb{R}^{1+n})
 } \nonumber\\
& \ \ \ \ \times  \prod^{\kappa_2+1}_{i=2} \|\Box_{k^{(i)}} u
 \|_{L^{\kappa_2}_{x_1} L^\infty_{x_2,...,x_n}L^\infty_t(\mathbb{R}^{1+n})
 }.
\label{est-int-2bb}
\end{align}
In view of the inclusion $L^{\kappa}_{x_1}
L^\infty_{x_2,...,x_n}L^\infty_t \bigcap L^{\infty}_{x,t}
 \subset  L^{\kappa_2}_{x_1}
L^\infty_{x_2,...,x_n}L^\infty_t$, we immediately have
\begin{align}
IV_1 & \lesssim  \rho_1(u) (\rho_2(u)+ \rho_3(u))^{\kappa_2}.
\label{est-int-2bc}
\end{align}
Noticing the fact that ${\rm supp} \psi_2 \subset \{\xi: \ |\xi_2|
\ge 2 |\xi_1|\}$ and applying \eqref{sm-int-1a}, we have
\begin{align}
IV_2 & \lesssim   \sum_{k \in \mathbb{Z}^n, \ |k_2| >4} \langle
k_2\rangle^{3/2} \sum_{\mathbb{S}^{(2)}_{2,2}} \|\Box_k
\left(\Box_{k^{(1)}} u ... \Box_{k^{(\kappa_2+1)}} u
\right)\|_{L^{(2+\kappa)/(1+\kappa)}_{t,x}(\mathbb{R}^{1+n})
}\nonumber\\
& \lesssim   \sum_{k \in \mathbb{Z}^n, \ |k_2| >4} \
\sum_{\mathbb{S}^{(2)}_{2,2}} \|\Box_k \left(\Box_{k^{(1)}} u ...
\Box_{k^{(\kappa_2+1)}} u
\right)\|_{L^{(2+\kappa)/(1+\kappa)}_{t,x}(\mathbb{R}^{1+n})
}\nonumber\\
& \lesssim  \rho_3(u)^{1+\kappa_2} . \label{est-int-2c}
\end{align}
The other terms in $\rho_1(\cdot)$ can be bounded in a similar way.
So, we have shown that
\begin{align}
\rho_1 \left(\mathscr{A} (\sum^n_{i=1} \lambda_i \partial_{x_i}
u^{\kappa_i+1}) \right) & \lesssim  \sum^n_{i=1} \left(\rho_1(u)
(\rho_2(u)+\rho_3(u))^{\kappa_i} + \rho_3(u)^{1+\kappa_i}\right) .
\label{est-int-2d}
\end{align}
We estimate $\rho_2 (\cdot)$. Denote
\begin{align}
\|u\|_{Z_i} = \sum_{k\in \mathbb{Z}^n} \langle
k\rangle^{1/2-1/\kappa} \|\Box_k u\|_{L^\kappa_{x_i}
L^\infty_{(x_j)_{j\not=i}}L^\infty_t(\mathbb{R}^{1+n})}.
\label{norm-Z}
\end{align}
We have
\begin{align}
\rho_2 \left(\mathscr{A} (\sum^n_{j=1} \lambda_i \partial_{x_i}
u^{\kappa_i+1}) \right) & \lesssim  \sum^n_{i=1}  \left
\|\mathscr{A} (\sum^n_{i=1} \lambda_i \partial_{x_i} u^{\kappa_i+1})
\right\|_{Z_j}. \label{est-rho-2a}
\end{align}
Due to the symmetry of $Z_1,...,Z_n$, it suffices to consider the
estimate of $\|\cdot\|_{Z_1}$. Recall that $k_{\max}=
|k_1|\vee...\vee |k_n|$. We have
\begin{align}
\|v\|_{Z_1} & \le \left(\sum_{k\in \mathbb{Z}^n, \, k_{\max}>4} +
\sum_{k\in \mathbb{Z}^n, \, k_{\max} \le 4} \right)\langle
k\rangle^{1/2-1/\kappa} \|\Box_k v\|_{L^\kappa_{x_1}
L^\infty_{x_2,...,x_n}L^\infty_t(\mathbb{R}^{1+n})} \nonumber\\
&:= \Gamma_1(v) + \Gamma_2(v). \label{est-rho-2b}
\end{align}
In view of Lemma \ref{sm-ef-int3} and H\"older's inequality,
\begin{align}
 \Gamma_2\left(\mathscr{A} \Big(\sum^n_{i=1} \lambda_i \partial_{x_i}
u^{\kappa_i+1} \Big) \right)  & \le \sum_{k\in \mathbb{Z}^n, \,
k_{\max} \le 4}
  \left\|\Box_k
\mathscr{A} \Big(\sum^n_{i=1} \lambda_i \partial_{x_i}
u^{\kappa_i+1} \Big) \right\|_{L^\kappa_{x_1}
L^\infty_{x_2,...,x_n}L^\infty_t(\mathbb{R}^{1+n})} \nonumber\\
& \lesssim  \sum^n_{i=1}
 \sum_{k^{(1)},...,k^{(\kappa_i+1)} \in \mathbb{Z}^n}
\left\|\Box_{k^{(1)}}u ... \Box_{k^{(\kappa_i+ 1)}}u
\right\|_{L^{\frac{2+\kappa_i}{1+\kappa_i}}_{t,x}(\mathbb{R}^{1+n})}
\nonumber\\
& \lesssim  \sum^n_{i=1}
 \sum_{k^{(1)},...,k^{(\kappa_i+1)} \in \mathbb{Z}^n}
\!\!\!\!\!\!\!\!\|\Box_{k^{(1)}}u\|_{L^{2+\kappa_i}_{t,x}(\mathbb{R}^{1+n})}
... \|\Box_{k^{(\kappa_i+ 1)}}u
\|_{L^{2+\kappa_i}_{t,x}(\mathbb{R}^{1+n})} \nonumber\\
& \lesssim  \sum^n_{i=1}  \rho_3(u)^{\kappa_i+1}. \label{est-rho-2c}
\end{align}
It is easy to see that
\begin{align}
\Gamma_1 (v)  & \le \left(\sum_{k\in \mathbb{Z}^n, \, |k_1|=
k_{\max}>4} +...+  \sum_{k\in \mathbb{Z}^n, \, |k_n| = k_{\max} > 4}
\right)\langle k\rangle^{1/2-1/\kappa} \|\Box_k v\|_{L^\kappa_{x_1}
L^\infty_{x_2,...,x_n}L^\infty_t(\mathbb{R}^{1+n})} \nonumber\\
&:= \Gamma^1_1(v) +...+ \Gamma^n_1(v) . \label{est-rho-2d}
\end{align}
Using \eqref{decomp}, Lemmas \ref{sm-ef-int3} and \ref{sm-ef-int4},
we have
\begin{align}
& \Gamma_1^1 \left(\mathscr{A} \Big(\sum^n_{i=1} \lambda_i
\partial_{x_i}
u^{\kappa_i+1} \Big) \right) \nonumber\\
& \lesssim \sum^n_{i=1} \sum_{k \in \mathbb{Z}^n, \ |k_1|
>4} \!\!\!\!\!\! \langle k_1 \rangle  \sum_{\mathbb{S}^{(1)}_{i,1}}
\|\Box_k \left(\Box_{k^{(1)}} u ... \Box_{k^{(\kappa_i+1)}} u
\right) \|_{L^{1}_{x_1}
L^2_{x_2,...,x_n}L^2_t(\mathbb{R}^{1+n})}  \nonumber\\
& \ \  + \sum^n_{i=1} \sum_{k \in \mathbb{Z}^n, \ |k_1| >4} \!\!\!
\langle k_1\rangle^{3/2} \sum_{\mathbb{S}^{(1)}_{i,2}} \| \Box_k
\left(\Box_{k^{(1)}} u ... \Box_{k^{(\kappa_i+1)}} u
\right)\|_{L^{(2+\kappa_i)/(1+\kappa_i)}_{t,x}(\mathbb{R}^{1+n}) }.
\label{nonhomo-2}
\end{align}
Using the same way as in \eqref{est-I} and \eqref{non-est-2}, one
easily sees that
\begin{align}
\Gamma_1^1 \left(\mathscr{A} \Big(\sum^n_{i=1} \lambda_i
\partial_{x_i}
u^{\kappa_i+1} \Big) \right) & \lesssim \sum^n_{i=1}\big(\rho_1(u)
(\rho_2(u)+ \rho_3(u))^{\kappa_i} + \rho_3(u)^{1+\kappa_i} \big).
\label{est-lam2}
\end{align}
We estimate $\Gamma^2_1(\cdot)$. By Lemmas \ref{sm-ef-int3} and
\ref{sm-ef-int4},
\begin{align}
& \Gamma_1^2 \left(\mathscr{A} \Big(\sum^n_{i=1} \lambda_i
\partial_{x_i}
u^{\kappa_i+1} \Big) \right) \nonumber\\
& \lesssim \sum^n_{i=1} \sum_{k \in \mathbb{Z}^n, \ |k_2|=k_{\max}
>4} \!\!\!\!\!\! \langle k \rangle^{1/2-1/\kappa}
\|\Box_k \left( \mathscr{A}
\partial_{x_i}
u^{\kappa_i+1}\right) \|_{L^{\kappa}_{x_1}
L^\infty_{x_2,...,x_n}L^\infty_t(\mathbb{R}^{1+n})}  \nonumber\\
& \lesssim \sum^n_{i=1} \sum_{k \in \mathbb{Z}^n, \ |k_2|=k_{\max}
>4} \!\!\!\!\!\! \langle k_2 \rangle  \sum_{\mathbb{S}^{(2)}_{i,1}}
\|\Box_k \left(\Box_{k^{(1)}} u ... \Box_{k^{(\kappa_i+1)}} u
\right) \|_{L^{1}_{x_2}
L^2_{x_1,x_3,...,x_n}L^2_t(\mathbb{R}^{1+n})}  \nonumber\\
& \ \  + \sum^n_{i=1} \sum_{k \in \mathbb{Z}^n, \ |k_2| >4} \!\!\!
\langle k_2\rangle^{3/2} \sum_{\mathbb{S}^{(2)}_{i,2}} \| \Box_k
\left(\Box_{k^{(1)}} u ... \Box_{k^{(\kappa_i+1)}} u
\right)\|_{L^{(2+\kappa_i)/(1+\kappa_i)}_{t,x}(\mathbb{R}^{1+n}) } .
\label{nonhomo-2a}
\end{align}
This reduces the same estimate as $\Gamma^1_1(\cdot)$. We easily see
that $\Gamma^i_1(\cdot)$ for $3\le i \le n$ can be controlled in a
similar way as $\Gamma^2_1(\cdot)$. Hence, we have shown that
\begin{align}
\left\|\mathscr{A} \Big(\sum^n_{i=1} \lambda_i \partial_{x_i}
u^{\kappa_i+1}\Big) \right\|_{Z_1} & \lesssim  \sum^n_{i=1}
\left(\rho_1(u) (\rho_2(u)+\rho_3(u))^{\kappa_i} +
\rho_3(u)^{1+\kappa_i}\right) . \label{est-Z-8}
\end{align}

For the estimates of $\rho_3(\mathscr{A} \partial_{x_i}
u^{\kappa_i+1})$, we have from \eqref{st-sm-mo-2} and Lemma
\ref{discret-deriv} that
\begin{align}
 \left \|\Box_k \mathscr{A} \partial_{x_i} f \right\|_{L^\infty_t L^2_x \, \cap \, L^{2+\kappa}_{t,x} (\mathbb{R}^{1+n}) } &
\lesssim    \|\Box_k \partial_{x_i} f\|_{
L^{(2+\kappa)/(1+\kappa)}_{t,x} (\mathbb{R}^{1+n}) } \nonumber\\
&  \lesssim   \langle k_i\rangle \|\Box_k  f\|_{
L^{(2+\kappa)/(1+\kappa)}_{t,x} (\mathbb{R}^{1+n}) }.
\label{st-sm-mo-2a}
\end{align}
Hence, using \eqref{decomp}, \eqref{st-sm-m-c9} and
\eqref{st-sm-m-c4}, we obtain that can be controlled by the right
hand side of \eqref{nonhomo-2}.
\begin{align}
 \rho_3(\mathscr{A}
\partial_{x_1} u^{\kappa_1+1}) & \lesssim \sum_{k \in \mathbb{Z}^n, \ |k_1| \le 4} \langle
k_1\rangle^{3/2} \sum_{k^{(1)},...,k^{(\kappa_1+1)} \in \mathbb{Z}^n
} \| \Box_k \left(\Box_{k^{(1)}} u ... \Box_{k^{(\kappa_1+1)}} u
\right)\|_{L^{(2+\kappa)/(1+\kappa)}_{t,x}(\mathbb{R}^{1+n}) }
\nonumber\\
& \ \ + \sum_{k \in \mathbb{Z}^n, \ |k_1| >4} \langle k_1\rangle
\sum_{\mathbb{S}^{(1)}_{1,1}} \|\Box_k \left(\Box_{k^{(1)}} u ...
\Box_{k^{(\kappa_1+1)}} u \right) \|_{L^{1}_{x_1}
L^2_{x_2,...,x_n}L^2_t(\mathbb{R}^{1+n})}  \nonumber\\
& \ \  +  \sum_{k \in \mathbb{Z}^n, \ |k_1| >4} \langle
k_1\rangle^{3/2} \sum_{\mathbb{S}^{(1)}_{1,2}} \| \Box_k
\left(\Box_{k^{(1)}} u ... \Box_{k^{(\kappa_1+1)}} u
\right)\|_{L^{(2+\kappa)/(1+\kappa)}_{t,x}(\mathbb{R}^{1+n}) }.
\label{nonh-1a}
\end{align}
By \eqref{est-I} and \eqref{non-est-2}, we have
\begin{align}
\rho_3(\mathscr{A} \partial_{x_1} u^{\kappa_1+1})
 & \lesssim
  \sum^n_{i=1} \left(\rho_1(u) (\rho_2(u)+\rho_3(u))^{\kappa_i} +
\rho_3(u)^{1+\kappa_i}\right). \label{est-rho34}
\end{align}
Hence, we have shown that
\begin{align}
\|\mathscr{T} u\|_X \lesssim \|u_0\|_{M^{1/2}_{2,1}}+ \sum^n_{i=1}
\|u\|^{1+\kappa_i}_X. \label{est-rho34}
\end{align}
Using a standard contraction mapping argument, we can finish the
proof of Theorem \ref{DNLS1-mod}. $\hfill \Box$

\section{Proof of Theorem \ref{DNLS-mod}} \label{pf-thm2}

Roughly speaking, we will prove our Theorem  \ref{DNLS-mod} by
following some ideas as in the proof of Theorem \ref{DNLS1-mod}.
However, due to the nonlinearity contains $u^{\kappa+1}$, and
$(\nabla u)^\nu$ and $u^\kappa (\nabla u)^\nu$ as special cases, the
proof of Theorem \ref{DNLS1-mod} can not be directly applied.  We
construct the space $X$ as follows. Denote
\begin{align*}
& \varrho^{(i)}_1(u) = \sum_{k\in \mathbb{Z}^n, \ |k_i|>4} \langle
k_i\rangle  \|\Box_k u\|_{L^{\infty}_{x_i}
L^2_{(x_j)_{j\not= i}}L^2_t(\mathbb{R}^{1+n})},\\
& \varrho^{(i)}_2(u) = \sum_{k\in \mathbb{Z}^n} \langle
k\rangle^{1/2-1/m} \|\Box_k u\|_{L^{m}_{x_i}
L^\infty_{(x_j)_{j\not= i}}L^\infty_t(\mathbb{R}^{1+n})},\\
& \varrho^{(i)}_3(u) = \sum_{k\in \mathbb{Z}^n} \langle k
\rangle^{1/2} \|\Box_k u\|_{ L^{2+m}_{x, t} \, \cap \,
L^\infty_tL^2_x (\mathbb{R}^{1+n})}.
\end{align*}
Put
\begin{align*}
X:= \left\{u\in \mathscr{S}'(\mathbb{R}^{1+n}): \ \|u\|_X:=
\sum^3_{\ell=1} \sum_{\alpha=0,1} \sum^n_{i,j=1}
\varrho^{(i)}_\ell(\partial^\alpha_{x_j} u) \le \delta \right\}.
\end{align*}
Considering the following mapping:
\begin{align*}
\mathscr{T}: u(t) \to S(t) u_0 - {\rm i} \mathscr{A} F(u,\bar{u},
\nabla u, \nabla\bar{u}),
\end{align*}
we will show that $\mathscr{T}: X\to X$ is a contraction mapping.

Since $\|u\|_X= \|\bar{u}\|_X$, we may assume, without loss of
generality that
\begin{align*}
 F(u,\bar{u},
\nabla u, \nabla\bar{u}) = F(u, \nabla u):= \sum_{m+1\le
\kappa+|\nu| <\infty} c_{\kappa \nu} u^\kappa (\nabla u)^{\nu},
\end{align*}
where $(\nabla u)^\nu = u^{\nu_1}_{x_1}...u^{\nu_n}_{x_n} $. For the
sake of convenience, we denote
\begin{align*}
v_1=...=v_\kappa=u, \ \ v_{\kappa+1}=...=v_{\kappa+\nu_1}= u_{x_1},
..., v_{\kappa+|\nu|-\nu_{n}+1}=...=v_{\kappa+|\nu|}= u_{x_n}.
\end{align*}
By \eqref{smo-eff.21}, for $\alpha=0,1$,
\begin{align*}
& \varrho^{(i)}_1(\partial^\alpha_{x_j} S(t)u_0) \lesssim \sum_{k\in
\mathbb{Z}^n, \ |k_i|>4} \langle k_i\rangle^{1/2} \langle k_j\rangle
\| \Box_k u_0\|_{L^{2} (\mathbb{R}^{n})} \le \|u_0\|_{
M^{3/2}_{2,1}}.
\end{align*}
By \eqref{st-sm-m-c3}, \eqref{st-sm-m-c1}, we have for $\alpha=0,1$,
\begin{align*}
& \varrho^{(i)}_2(\partial^\alpha_{x_j} S(t)u_0) +
\varrho^{(i)}_3(\partial^\alpha_{x_j} S(t)u_0) \lesssim \|u_0\|_{
M^{3/2}_{2,1}}.
\end{align*}
Hence,
\begin{align*}
\|S(t)u_0\|_X   \lesssim \|u_0\|_{ M^{3/2}_{2,1}}.
\end{align*}
In order to estimate $\varrho^{(i)}_1(\mathscr{A}
\partial^\alpha_{x_j}(v_1... v_{\kappa+|\nu|}))$, $i,j=1,...,n$,  it suffices to estimate $\varrho^{(1)}_1(\mathscr{A} \partial^\alpha_{x_1}(v_1...
v_{\kappa+|\nu|}))$ and $\varrho^{(1)}_1(\mathscr{A}
\partial^\alpha_{x_2}(v_1... v_{\kappa+|\nu|}))$.  Similarly as in
\eqref{decomp}, we will use the decomposition
\begin{align}
\Box_k (v_1...v_{\kappa+|\nu|})
 = & \sum_{\mathbb{S}^{(i)}_1} \Box_k
\left(\Box_{k^{(1)}} v_1 ... \Box_{k^{(\kappa+|\nu|)}}
v_{\kappa+|\nu|}
\right) \nonumber\\
&  +  \sum_{\mathbb{S}^{(i)}_2} \Box_k \left(\Box_{k^{(1)}} v_1 ...
\Box_{k^{(\kappa+|\nu|)}} v_{\kappa+|\nu|} \right), \label{decomp2}
\end{align}
where
\begin{align*}
& \mathbb{S}^{(i)}_1:= \{(k^{(1)},..., k^{(\kappa+|\nu
|)}): \
|k^{(1)}_i|\vee...\vee |k^{(\kappa+|\nu|)}_i|>4\}, \\
& \mathbb{S}^{(i)}_2:= \{(k^{(1)},..., k^{(\kappa+|\nu|)}): \
|k^{(1)}_i|\vee...\vee |k^{(\kappa+|\nu|)}_i|\le 4\}.
\end{align*}
In view of \eqref{1-sm-mod} and \eqref{st-sm-mo-4},
\begin{align}
& \varrho^{(1)}_1(\mathscr{A} \partial^\alpha_{x_1}(v_1...
v_{\kappa+|\nu|})) \nonumber\\
 & \lesssim \sum_{k \in \mathbb{Z}^n, \
|k_1|
>4} \langle k_1\rangle \sum_{\mathbb{S}^{(1)}_1} \|\Box_k
\left(\Box_{k^{(1)}} v_1 ... \Box_{k^{(\kappa+|\nu|)}}
v_{\kappa+|\nu|} \right) \|_{L^{1}_{x_1}
L^2_{x_2,...,x_n}L^2_t(\mathbb{R}^{1+n})}  \nonumber\\
& \ \  +  \sum_{k \in \mathbb{Z}^n, \ |k_1| >4} \langle
k_1\rangle^{3/2} \sum_{\mathbb{S}^{(1)}_2} \| \Box_k
\left(\Box_{k^{(1)}} v_1 ... \Box_{k^{(\kappa+|\nu|)}}
v_{\kappa+|\nu|} \right)\|_{L^{\frac{\kappa+|\nu|+1}{\kappa+|\nu|}}_{t,x}(\mathbb{R}^{1+n}) } \nonumber\\
& := I +II . \label{gnonhomo-1}
\end{align}
Similar to \eqref{est-I},
\begin{align}
I & \lesssim    \sum_{k^{(1)}\in \mathbb{Z}^n, \ |k^{(1)}_1|>2}
\langle k^{(1)}_1\rangle  \|\Box_{k^{(1)}}
v_1\|_{L^{\infty}_{x_1} L^2_{x_2,...,x_n}L^2_t(\mathbb{R}^{1+n})} \nonumber\\
& \ \ \times  \sum_{k^{(2)},...,k^{(\kappa+|\nu|)}\in \mathbb{Z}^n}
\  \prod^{\kappa+|\nu|}_{i=2} \|\Box_{k^{(i)}} v_i
\|_{L^{\kappa+|\nu|-1}_{x_1}
L^\infty_{x_2,...,x_n}L^\infty_t(\mathbb{R}^{1+n})}. \label{gest-I}
\end{align}
By H\"older's inequality and Lemma \ref{qnorm:pnorm},
\begin{align}
&  \|\Box_{k^{(i)}} v_i \|_{L^{\kappa+|\nu|-1}_{x_1}
L^\infty_{x_2,...,x_n}L^\infty_t(\mathbb{R}^{1+n})} \nonumber\\
& \le  \|\Box_{k^{(i)}} v_i
\|^{\frac{m}{\kappa+|\nu|-1}}_{L^{m}_{x_1}
L^\infty_{x_2,...,x_n}L^\infty_t(\mathbb{R}^{1+n})} \|\Box_{k^{(i)}}
v_i \|^{1- \frac{m}{\kappa+|\nu|-1}}_{
L^\infty_{x,t}(\mathbb{R}^{1+n}) } \nonumber\\
& \lesssim  \|\Box_{k^{(i)}} v_i
\|^{\frac{m}{\kappa+|\nu|-1}}_{L^{m}_{x_1}
L^\infty_{x_2,...,x_n}L^\infty_t(\mathbb{R}^{1+n})} \|\Box_{k^{(i)}}
v_i \|^{1- \frac{m}{\kappa+|\nu|-1}}_{ L^\infty_t
L^2_{x}(\mathbb{R}^{1+n}) }. \label{gest-I-con}
\end{align}
Hence, noticing that $v_i=u$ or $v_i=u_{x_j}$,  we have from
\eqref{gest-I} and \eqref{gest-I-con},
\begin{align}
I \lesssim \|u\|_X^{\kappa+|\nu|}. \label{I-control}
\end{align}
Similar to \eqref{non-est-2}, we see that $|k_1| \le C$ in the
summation of $II$. Again, in view of H\"older's inequality and Lemma
\eqref{qnorm:pnorm},
\begin{align}
\|\Box_{k^{(1)}} v_1 ... \Box_{k^{(\kappa+|\nu|)}}
v_{\kappa+|\nu|}\|_{L^{\frac{\kappa+|\nu|+1}{\kappa+|\nu|}}_{x,t}
(\mathbb{R}^{1+n})}  & \le \prod^{\kappa+|\nu|}_{i=1}
\|\Box_{k^{(i)}} v_i
\|_{L^{\kappa+|\nu|+1}_{x,t} (\mathbb{R}^{1+n})} \nonumber\\
& \lesssim \prod^{\kappa+|\nu|}_{i=1}  \|\Box_{k^{(i)}} v_i
\|_{L^{2+m}_{x,t} \, \bigcap \, L^\infty_t L^2_{x}
(\mathbb{R}^{1+n})}. \label{gest-II-con}
\end{align}
Hence, using a similar way as in \eqref{non-est-2},
\begin{align}
II \lesssim \|u\|_X^{\kappa+|\nu|}. \label{II-control}
\end{align}
We now give the estimate of $\varrho^{(1)}_1(\mathscr{A}
\partial^\alpha_{x_2}(v_1... v_{\kappa+|\nu|}))$.
Since we have obtained the estimate in the case $\alpha=0$, it suffices to consider the case $\alpha=1$.
 Let $\psi_i$ $(i=1,2)$ be as in Lemma
\ref{sm-ef-int2} and $P_i= \mathscr{F}^{-1} \psi_i \mathscr{F}$. We
have
\begin{align}
& \varrho^{(1)}_1(\mathscr{A}
\partial_{x_2}(v_1... v_{\kappa+|\nu|})) \nonumber\\
& \le  \sum_{k\in \mathbb{Z}^n, \ |k_1|>4} \langle k_1\rangle  \|P_1
\Box_k (\mathscr{A}
\partial_{x_2}(v_1... v_{\kappa+|\nu|}))\|_{L^\infty_{x_1} L^2_{x_2,...,x_n}L^2_t} \nonumber\\
& \ \ + \sum_{k\in \mathbb{Z}^n, \ |k_1|>4} \langle k_1\rangle
\|P_2 \Box_k (\mathscr{A}
\partial_{x_2}(v_1... v_{\kappa+|\nu|}))\|_{L^\infty_{x_1} L^2_{x_2,...,x_n}L^2_t} \nonumber\\
& := III+IV. \label{rho112-1}
\end{align}
Using the decomposition \eqref{decomp2},
\begin{align}
III & \le  \sum_{k\in \mathbb{Z}^n, \ |k_1|>4} \langle k_1\rangle
\sum_{\mathbb{S}^{(1)}_1}   \|P_1 \Box_k (\mathscr{A}
\partial_{x_2}(\Box_{k^{(1)}} v_1 ... \Box_{k^{(\kappa+|\nu|)}}
v_{\kappa+|\nu|} ))\|_{L^\infty_{x_1} L^2_{x_2,...,x_n}L^2_t} \nonumber\\
& \ \ + \sum_{k\in \mathbb{Z}^n, \ |k_1|>4} \langle k_1\rangle
\sum_{\mathbb{S}^{(1)}_2}   \|P_1 \Box_k (\mathscr{A}
\partial_{x_2}(\Box_{k^{(1)}} v_1 ... \Box_{k^{(\kappa+|\nu|)}}
v_{\kappa+|\nu|}))\|_{L^\infty_{x_1} L^2_{x_2,...,x_n}L^2_t} \nonumber\\
& := III_1 + III_2 \label{rho112-2}.
\end{align}
By Lemma \ref{sm-ef-int2},
\begin{align}
III_1 & \lesssim   \sum_{\mathbb{S}^{(1)}_1} \sum_{k\in
\mathbb{Z}^n, \ |k_1|>4} \langle k_1\rangle     \| \Box_k
(\Box_{k^{(1)}} v_1 ... \Box_{k^{(\kappa+|\nu|)}} v_{\kappa+|\nu|}
)\|_{L^1_{x_1} L^2_{x_2,...,x_n}L^2_t}.
 \label{rho112-3}
\end{align}
By symmetry, we may assume $|k^{(1)}_1|= \max (|k^{(1)}_1|,...,
|k^{(\kappa+|\nu|)}_1|)$ in $\mathbb{S}^{(1)}_1$. Hence,
\begin{align}
III_1 & \lesssim   \sum_{\mathbb{S}^{(1)}_1, \ |k^{(1)}_1|>4}
\langle k^{(1)}_1\rangle   \|\Box_{k^{(1)}} v_1\|_{L^\infty_{x_1}
L^2_{x_2,...,x_n}L^2_t} \prod^{\kappa+|\nu|}_{i=2} \| \Box_{k^{(i)}}
v_i \|_{L^{\kappa+|\nu|-1}_{x_1}
L^\infty_{x_2,...,x_n}L^\infty_t} \nonumber\\
& \lesssim   \varrho^{(1)}_1 (v_1) \prod^{\kappa+|\nu|}_{i=2}
(\varrho^{(1)}_2 (v_i) + \varrho^{(1)}_3 (v_i)) \lesssim
\|u\|^{\kappa+|\nu|}_X .
 \label{rho112-4}
\end{align}
Applying \eqref{sm-int-1a} and using a similar way as in
\eqref{est-int-2},
\begin{align}
III_2 & \lesssim  \!\!\! \sum_{k \in \mathbb{Z}^n, \ |k_1| >4, \
|k_2| \lesssim |k_1|} \!\!\! \langle k_1\rangle^{3/2}
\sum_{\mathbb{S}^{(1)}_2} \|\Box_k (\Box_{k^{(1)}} v_1 ...
\Box_{k^{(\kappa+|\nu|)}} v_{\kappa+|\nu|}
)\|_{L^{(2+m)/(1+m)}_{t,x}(\mathbb{R}^{1+n})
}\nonumber\\
& \lesssim  \prod^{\kappa+|\nu|}_{i=1} \varrho^{(1)}_3(v_i) \le
\|u\|^{\kappa+|\nu|}_X. \label{rho112-5}
\end{align}
So, we have shown that
\begin{align}
III & \lesssim  \|u\|^{\kappa+|\nu|}_X. \label{rho112-6}
\end{align}
Now we estimate $IV$. Using the decomposition \eqref{decomp2},
\begin{align}
IV & \le  \sum_{k\in \mathbb{Z}^n, \ |k_1|>4} \langle k_1\rangle
\sum_{\mathbb{S}^{(2)}_1}   \|P_2 \Box_k (\mathscr{A}
\partial_{x_2}(\Box_{k^{(1)}} v_1 ... \Box_{k^{(\kappa+|\nu|)}}
v_{\kappa+|\nu|} ))\|_{L^\infty_{x_1} L^2_{x_2,...,x_n}L^2_t} \nonumber\\
& \ \ + \sum_{k\in \mathbb{Z}^n, \ |k_1|>4} \langle k_1\rangle
\sum_{\mathbb{S}^{(2)}_2}   \|P_2 \Box_k (\mathscr{A}
\partial_{x_2}(\Box_{k^{(1)}} v_1 ... \Box_{k^{(\kappa+|\nu|)}}
v_{\kappa+|\nu|}))\|_{L^\infty_{x_1} L^2_{x_2,...,x_n}L^2_t} \nonumber\\
& := IV_1 + IV_2 \label{rho112-7}.
\end{align}
By Lemma \ref{sm-ef-int2},
\begin{align}
IV_1 & \lesssim   \sum_{\mathbb{S}^{(2)}_1} \sum_{k\in \mathbb{Z}^n,
\ |k_2|>4} \langle k_2\rangle     \| \Box_k (\Box_{k^{(1)}} v_1 ...
\Box_{k^{(\kappa+|\nu|)}} v_{\kappa+|\nu|} )\|_{L^1_{x_1}
L^2_{x_2,...,x_n}L^2_t}.
 \label{rho112-8}
\end{align}
In view of the symmetry, one can bound $IV_1$ by using the same way
as that of $III_1$ and as in
\eqref{est-int-2b}--\eqref{est-int-2bc}:
\begin{align}
IV_1 & \lesssim  \|u\|^{\kappa+|\nu|}_X. \label{rho112-9}
\end{align}
For the estimate of $IV_2$, we apply \eqref{sm-int-1a},
\begin{align}
IV_2 & \lesssim  \!\!\! \sum_{k \in \mathbb{Z}^n, \ |k_1| >4} \!\!\!
\langle k_1\rangle^{1/2} \langle k_2\rangle
\sum_{\mathbb{S}^{(2)}_2} \|P_2 \Box_k (\Box_{k^{(1)}} v_1 ...
\Box_{k^{(\kappa+|\nu|)}} v_{\kappa+|\nu|}
)\|_{L^{(2+m)/(1+m)}_{t,x}(\mathbb{R}^{1+n})
}\nonumber\\
& \lesssim  \sum_{\mathbb{S}^{(2)}_2} \|\Box_{k^{(1)}} v_1 ...
\Box_{k^{(\kappa+|\nu|)}} v_{\kappa+|\nu|}
\|_{L^{(2+m)/(1+m)}_{t,x}(\mathbb{R}^{1+n})}\lesssim
\|u\|^{\kappa+|\nu|}_X. \label{rho112-10}
\end{align}
Hence, in view of \eqref{rho112-9} and  \eqref{rho112-10}, we have
\begin{align}
IV \lesssim \|u\|^{\kappa+|\nu|}_X. \label{rho112-11}
\end{align}
Collecting \eqref{I-control}, \eqref{II-control}, \eqref{rho112-6},
\eqref{rho112-11}, we have shown that
\begin{align}
\sum_{\alpha=0,1} \sum^n_{i,j=1}\varrho_1^{(i)}
(\mathscr{A}\partial^\alpha_{x_j} (u^\kappa (\nabla u)^\nu))
\lesssim \|u\|_X^{\kappa+|\nu|}. \label{lambda1-control}
\end{align}

\begin{lem}\label{lem5.1}
Let $s\ge 0$,  $1\le p, p_i, \gamma, \gamma_i \le \infty$ satisfy
\begin{align}
\frac{1}{p}= \frac{1}{p_1}+...+\frac{1}{p_N}, \quad
\frac{1}{\gamma}= \frac{1}{\gamma_1}+...+   \frac{1}{\gamma_N}.
\label{p-gamma}
\end{align}
Then
\begin{align}
\sum_{k\in \mathbb{Z}^n} \langle k_1\rangle^{s} \left\|\Box_k (u_1
... u_N) \right\|_{L^{\gamma}_t L^p_x (\mathbb{R}^{1+n})} & \lesssim
\prod^N_{i=1} \left(\sum_{k\in \mathbb{Z}^n} \langle
k_1\rangle^{s}\|\Box_k u_i\|_{L^{\gamma_i}_t L^{p_i}_x
(\mathbb{R}^{1+n})}\right). \label{P7}
\end{align}
\end{lem}
{\bf Proof.} See \cite{WaHe}, Lemma 7.1.  $\hfill\Box$

\medskip

Next, we consider the estimates of $\varrho_2^{(1)} (\mathscr{A}
(u^\kappa (\nabla u)^\nu))$ and $\varrho_3^{(1)} (\mathscr{A}
(u^\kappa (\nabla u)^\nu))$. In view of \eqref{st-sm-m-c9} and
\eqref{st-sm-mo-6},
\begin{align}
\sum_{j=2,3} \varrho_j^{(1)} (\mathscr{A} (u^\kappa (\nabla u)^\nu))
\lesssim  \sum_{k\in \mathbb{Z}^n} \langle k \rangle^{1/2} \|\Box_k
(u^\kappa (\nabla
u)^\nu)\|_{L^{\frac{2+m}{1+m}}_{t,x}(\mathbb{R}^{1+n}) }.
\label{lambda23-control}
\end{align}
We use Lemma \ref{lem5.1} to control the right hand side of
\eqref{lambda23-control}:
\begin{align}
& \sum_{k\in \mathbb{Z}^n} \langle k\rangle^{1/2} \|\Box_k
(v_1...v_{\kappa+|\nu|})
\|_{L^{\frac{2+m}{1+m}}_{t,x}(\mathbb{R}^{1+n})} \nonumber\\
& \lesssim \prod^{m+1}_{i=1} \left(\sum_{k\in \mathbb{Z}^n} \langle
k\rangle^{1/2}\|\Box_k v_i\|_{L^{2+m}_{t,x}
(\mathbb{R}^{1+n})}\right)  \prod^{\kappa+|\nu|}_{i=m+2}
\left(\sum_{k\in \mathbb{Z}^n} \langle k \rangle^{1/2}\|\Box_k
v_i\|_{L^{\infty}_{t,x} (\mathbb{R}^{1+n})}\right) \nonumber\\
& \lesssim \prod^{m+1}_{i=1} \left(\sum_{k\in \mathbb{Z}^n} \langle
k \rangle^{1/2}\|\Box_k v_i\|_{L^{2+m}_{t,x}
(\mathbb{R}^{1+n})}\right) \prod^{\kappa+|\nu|}_{i=m+2}
\left(\sum_{k\in \mathbb{Z}^n} \langle k \rangle^{1/2}\|\Box_k
v_i\|_{L^{\infty}_{t}L^2_{x} (\mathbb{R}^{1+n})}\right) \nonumber\\
& \lesssim \prod^{\kappa+|\nu|}_{i=1} \varrho^{(1)}_3(v_i) \le
\|u\|^{\kappa+|\nu|}_X.  \label{lambda23-cont}
\end{align}
We estimate $\varrho_2^{(1)} (\mathscr{A}
\partial_{x_1}(u^\kappa (\nabla u)^\nu))$. Recall that $k_{\max}= |k_1|\vee ...\vee
|k_n|$.
\begin{align}
&  \varrho^{(1)}_2 (\mathscr{A} \partial_{x_1}(v_1...
v_{\kappa+|\nu|})) \nonumber\\
& \lesssim \sum_{k \in \mathbb{Z}^n, \
k_{\max}
>4} \langle
k \rangle^{1/2-1/m}  \|\Box_k \mathscr{A} \partial_{x_1} \left( v_1
... v_{\kappa+|\nu|} \right) \|_{L^{m}_{x_1}
L^\infty_{x_2,...,x_n}L^\infty_t(\mathbb{R}^{1+n})}  \nonumber\\
& \ \ \ \ +  \sum_{k \in \mathbb{Z}^n, \ k_{\max} \le 4} \langle k
\rangle^{1/2-1/m} \|\Box_k \mathscr{A} \partial_{x_1} \left( v_1 ...
v_{\kappa+|\nu|} \right) \|_{L^{m}_{x_1}
L^\infty_{x_2,...,x_n}L^\infty_t(\mathbb{R}^{1+n})} \nonumber\\
& := V+VI . \label{gnho-final}
\end{align}
By \eqref{st-sm-mo-6} and Lemma \ref{lem5.1}, we have
\begin{align}
VI \lesssim \sum_{k\in \mathbb{Z}^n}  \|\Box_k
(v_1...v_{\kappa+|\nu|})
\|_{L^{\frac{2+m}{1+m}}_{t,x}(\mathbb{R}^{1+n})} \lesssim
\|u\|^{\kappa+|\nu|}_X. \label{gnonhomo-B}
\end{align}
It is easy to see that
\begin{align}
V & \lesssim \left(\sum_{k \in \mathbb{Z}^n, \ |k_1|= k_{\max}
>4}+...+ \sum_{k \in \mathbb{Z}^n, \ |k_n|= k_{\max}>4} \right)
\langle k \rangle^{1/2-1/m} \nonumber\\
& \ \ \ \ \ \times \|\Box_k \mathscr{A}
\partial_{x_1} \left( v_1 ... v_{\kappa+|\nu|} \right)
\|_{L^{m}_{x_1}
L^\infty_{x_2,...,x_n}L^\infty_t(\mathbb{R}^{1+n})}:=
\Upsilon_1(u)+...+ \Upsilon_n(u). \label{gnho-finala}
\end{align}
Applying the decomposition \eqref{decomp2} and Lemmas
\ref{sm-ef-int3} and \ref{sm-ef-int4}, we obtain that
\begin{align}
\Upsilon_1(u) & \lesssim \sum_{k \in \mathbb{Z}^n, \ |k_1|
>4} \langle k_1\rangle  \sum_{\mathbb{S}^{(1)}_1} \|\Box_k
\left(\Box_{k^{(1)}} v_1 ... \Box_{k^{(\kappa+|\nu|)}}
v_{\kappa+|\nu|} \right) \|_{L^{1}_{x_1}
L^2_{x_2,...,x_n}L^2_t(\mathbb{R}^{1+n})}  \nonumber\\
& \ \  +  \sum_{k \in \mathbb{Z}^n, \ |k_1| >4} \langle
k_1\rangle^{3/2} \sum_{\mathbb{S}^{(1)}_2} \| \Box_k
\left(\Box_{k^{(1)}} v_1 ... \Box_{k^{(\kappa+|\nu|)}}
v_{\kappa+|\nu|}
\right)\|_{L^{\frac{\kappa+|\nu|+1}{\kappa+|\nu|}}_{t,x}(\mathbb{R}^{1+n})
}, \label{gnonhomo-A}
\end{align}
which reduces to the case $\alpha=1$ in \eqref{gnonhomo-1}. So,
\begin{align}
\Upsilon_1(u)  & \lesssim  \|u\|^{\kappa+|\nu|}_X.
\label{gnonh-A+Ba}
\end{align}
Again, in view of Lemmas \ref{sm-ef-int3} and \ref{sm-ef-int4},
\begin{align}
\Upsilon_2(u) & \lesssim \sum_{k \in \mathbb{Z}^n, \ |k_2|
>4} \langle k_2\rangle  \sum_{\mathbb{S}^{(2)}_1} \|\Box_k
\left(\Box_{k^{(1)}} v_1 ... \Box_{k^{(\kappa+|\nu|)}}
v_{\kappa+|\nu|} \right) \|_{L^{1}_{x_2}
L^2_{x_1,x_3,...,x_n}L^2_t(\mathbb{R}^{1+n})}  \nonumber\\
& \ \  +  \sum_{k \in \mathbb{Z}^n, \ |k_2| >4} \langle
k_2\rangle^{3/2} \sum_{\mathbb{S}^{(2)}_2} \| \Box_k
\left(\Box_{k^{(1)}} v_1 ... \Box_{k^{(\kappa+|\nu|)}}
v_{\kappa+|\nu|}
\right)\|_{L^{\frac{\kappa+|\nu|+1}{\kappa+|\nu|}}_{t,x}(\mathbb{R}^{1+n})
}, \label{gnonhomo-Aa}
\end{align}
which reduces to the same estimate as $\Upsilon_1(u) $.  Using the
same way as $\Upsilon_2(u)$, we can get the estimates of
$\Upsilon_3(u),...,\Upsilon_n(u)$. So,
\begin{align}
 \varrho^{(1)}_2 (\mathscr{A} \partial_{x_1}(v_1...
v_{\kappa+|\nu|})) & \lesssim  \|u\|^{\kappa+|\nu|}_X.
\label{gnonh-A+B}
\end{align}
We need to further bound $ \varrho^{(1)}_2 (\mathscr{A}
\partial_{x_i}(v_1... v_{\kappa+|\nu|}))$, $i=2,...,n$, which is essentially the same as $ \varrho^{(1)}_2 (\mathscr{A}
\partial_{x_1}(v_1... v_{\kappa+|\nu|}))$. Indeed, it is easy to see
that  \eqref{gnho-final} holds if we substitute $\partial_{x_1}$
with $\partial_{x_i}$. Moreover, using Lemmas \ref{lem5.1},
\ref{sm-ef-int3} and \ref{sm-ef-int4}, we easily get that
\begin{align}
\varrho^{(1)}_2 \left(\mathscr{A} (
\partial_{x_i} ( v_1... v_{\kappa+|\nu|})) \right) & \lesssim \|u\|_X^{\kappa+|\nu|}. \label{es-rho21-8}
\end{align}

By Lemma \ref{discret-deriv}, \eqref{st-sm-mo-2}, we see that
\begin{align}
\|\Box_k \mathscr{A} \partial_{x_1} f\|_{L^\infty_t L^2 \, \cap \,
 L^{2+m}_{t,x}(\mathbb{R}^{1+n})} \lesssim \langle k_1\rangle
 \|\Box_k f\|_{L^{\frac{2+m}{1+m}}_{t,x}(\mathbb{R}^{1+n})}.
\label{st-sm-m-7a}
\end{align}
Hence, in view of \eqref{st-sm-mo-7} and \eqref{st-sm-mo-3},
repeating the procedure as in the estimates of $\rho_3(u)$ in
Theorem \ref{DNLS1-mod}, $\varrho^{(1)}_3 (\mathscr{A}
\partial_{x_1}(v_1... v_{\kappa+|\nu|}))$ can be controlled by the
right hand side of \eqref{gnonhomo-A} and \eqref{gnonhomo-B}.
Summarizing the estimates as in the above, we have shown
that\footnote{Notice that $|c_\beta|\le C^{|\beta|}$.}
\begin{align}
\|\mathscr{T} u\|_X \le C \|u_0\|_{M^{3/2}} + \sum_{m+1\le \ell
<\infty } \ell^{2n+2} C^\ell \|u\|^{\ell}_X. \label{est-full}
\end{align}
Applying a standard contraction mapping argument, we can prove our
result.

\section{Proofs of Theorems \ref{DNLS1-modm2} and \ref{DNLS-modm2}}

{\bf Proof of Theorem \ref{DNLS1-modm2}.} For convenience, we denote
\begin{align}
 & \rho_1(u) =  \sum^n_{i=1} \sum_{k\in \mathbb{Z}^n, \ |k_i|>4}
\langle k_i\rangle^{2}  \left\| \Box_k u \right\|_{L^\infty_{x_i}
L^2_{(x_j)_{j\not=i}}L^2_t(\mathbb{R}^{1+n})} \nonumber\\
&  \rho_2(u) =  \sum^n_{i=1} \sum_{k\in \mathbb{Z}^n}  \left\|
\Box_k u \right\|_{L^{\kappa}_{x_i}
L^\infty_{(x_j)_{j\not=i}}L^\infty_t(\mathbb{R}^{1+n})} \nonumber\\
&  \rho_3(u) = \sum_{k\in \mathbb{Z}^n}  \langle k \rangle^{3/2}
\left\| \Box_k u \right\|_{L^\infty_t L^2_x  \bigcap L^3_t L^{6}_{x}
(\mathbb{R}^{1+n})}. \nonumber
\end{align}
Comparing the definitions of $\rho_i (u)$ with those of Section
\ref{pf-thm1}, we see that here we drop the regularity $\langle
k_i\rangle^{1/2-1/\kappa}$ in $\rho_2 (u)$ and we add $1$-order
regularity in $\rho_1(u)$ and $\rho_3(u)$. The estimates for
$\rho_1(\mathscr{T} u)$ and  $\rho_3(\mathscr{T} u)$ can be shown by
following the same way as in the Section \ref{pf-thm1} (It is worth
to notice that in Section \ref{pf-thm1}, when we estimate
$\rho_1(\mathscr{T} u)$ and  $\rho_3(\mathscr{T} u)$, we can replace
$\rho_2(u)$ defined here to substitute that in Section
\ref{pf-thm1}). We also need to point out that for $n\ge 2$, $2/3 <
n(1/2-1/6)$ and so, $\|\cdot \|_{L^3_t L^{6}_{x}
(\mathbb{R}^{1+n})}$ is a Strichartz norm. Moreover,
$$ \|\Box_k u\|_{L^{2+p}_{x,t}} \lesssim
\|\Box_k u\|_{L^\infty_t L^2_x  \bigcap L^3_t L^{6}_{x}
(\mathbb{R}^{1+n})}
$$
uniformly holds for all $k\in \mathbb{Z}^n$ and $2\le p\le \infty$.

Noticing that in the proof of Theorem \ref{DNLS1-mod}, we do not
know if the following two inequalities hold for $m=2$,
\begin{align}
& \left \|\Box_k \mathscr{A} \partial_{x_1}  f \right\|_{L^m_{x_1}
L^\infty_{x_2,...,x_n} L^\infty_t (\mathbb{R}^{1+n})}  \lesssim
\langle k_1\rangle^{1/2+1/m} \| \Box_k f\|_{L^1_{x_1}
L^2_{x_2,...,x_n} L^2_t (\mathbb{R}^{1+n})}, \label{st-sm-m-c6a}\\
& \left \|\Box_k \partial_{x_i}  \mathscr{A}  f \right\|_{L^m_{x_1}
L^\infty_{x_2,...,x_n} L^\infty_t (\mathbb{R}^{1+n})} \lesssim
\langle k_i\rangle^{1/2}  \langle k_1\rangle^{1/m} \|
 \Box_k  f\|_{L^1_{x_i}
L^2_{(x_j)_{j\not=i}}L^2_t (\mathbb{R}^{1+n})}. \label{sm-int-2a}
\end{align}
So, in the case $m=2$, we need to find another way to estimate
$\rho_2(\mathscr{T} u)$. Our solution is to apply the following
estimate as in \eqref{st-sm-mo-6}:
\begin{align}
& \left \|\Box_k \mathscr{A}  f \right\|_{L^2_{x_1}
L^\infty_{x_2,...,x_n} L^\infty_t (\mathbb{R}^{1+n})}  \lesssim
\langle k_1\rangle^{1/2} \| \Box_k f\|_{L^1_{t}
L^2_{x}(\mathbb{R}^{1+n})}. \label{st-sm-m-c6b}
\end{align}
It follows that for any $\kappa\ge 2$,
\begin{align}
& \sum_{k\in \mathbb{Z}^n}\left \|\Box_k \mathscr{A} \partial_{x_i}
u^{\kappa+1} \right\|_{L^2_{x_1} L^\infty_{x_2,...,x_n} L^\infty_t
(\mathbb{R}^{1+n})}  \lesssim \sum_{k\in \mathbb{Z}^n} \langle
k\rangle^{3/2} \| \Box_k u^{\kappa+1}\|_{L^1_{t}
L^2_{x}(\mathbb{R}^{1+n})}. \label{st-sm-m-c6c}
\end{align}
Using Lemma \ref{lem5.1}, one has that
\begin{align}
\sum_{k\in \mathbb{Z}^n}\left \|\Box_k \mathscr{A} \partial_{x_i}
u^{\kappa+1} \right\|_{L^2_{x_1} L^\infty_{x_2,...,x_n} L^\infty_t
(\mathbb{R}^{1+n})}  & \lesssim \left(\sum_{k\in \mathbb{Z}^n}
\langle k\rangle^{3/2} \| \Box_k u\|_{L^{\kappa+1}_{t}
L^{2(\kappa+1)}_{x}(\mathbb{R}^{1+n})} \right)^{\kappa+1}\nonumber\\
& \lesssim \left(\sum_{k\in \mathbb{Z}^n} \langle k\rangle^{3/2} \|
\Box_k u\|_{L^{3}_{t}
L^{6}_{x} \cap L^\infty_{x,t} (\mathbb{R}^{1+n})}\right)^{\kappa+1} \nonumber\\
& \lesssim  \rho_3(u)^{1+\kappa}. \label{st-sm-m-c6d}
\end{align}
Using \eqref{st-sm-m-c6d}, the estimates of $\rho_2(\mathscr{T} u)$
is also obtained. $\hfill \Box$\\

\noindent {\bf Proof of Theorem \ref{DNLS-modm2}.} We can follow the
proof of Theorems \ref{DNLS1-modm2} and \ref{DNLS-mod} to get the
proof and we omit the details of the proof. $\hfill \Box$

\begin{appendix}
\section{Appendix}

In this section, we generalize the Christ-Kiselev Lemma
\cite{Ch-Kis} to anisotropic Lebesgue spaces. Our idea follows
Molinet and Ribaud \cite{Mo-Ri}, and Smith and Sogge \cite{SS}.
Denote
\begin{align}
& T f(t)=\int_{-\infty}^{\infty}K(t, t')f(t')dt', \ \ \   T_{re}
f(t)=\int_{0}^{t} K(t, t')f(t')dt'.  \label{A1}
\end{align}
If $T: \ Y_1\to X_1$ implies that  $T_{re}: \ Y_1\to X_1$, then $T:
\ Y_1\to X_1$ is said to be a well restriction operator.

\begin{prop} \label{propA.1}
Let $T$ be as in \eqref{A1}. We have the following results.
\begin{itemize}
     \item[\rm (1)] If $\wedge^3_{i=1} p_i > (\vee^3_{i=1} q_i) \vee ( q_1q_3 /q_2)$,
     then $T: L_{x_1}^{q_1}L_{x_2}^{q_2}L_t^{q_3}(\mathbb{R}^{3})
     \to
     L_{x_1}^{p_1}L_{x_2}^{p_2}L_t^{p_3} (\mathbb{R}^{3})$
     is a well restriction operator.
\item[\rm (2)] If  $p_1 > (\vee^3_{i=1} q_i) \vee ( q_1q_3/q_2)$,
then $T: \ L_{x_1}^{q_1}L_{x_2}^{q_2}L_t^{q_3}(\mathbb{R}^{3}) \to
L_{t}^{p_1}L_{x_1}^{p_2}L_{x_2}^{p_3}(\mathbb{R}^{3})$ is a well
restriction operator.
\item[\rm (3)] If $q_1< \wedge^3_{i=1} p_i$, then
$T: L_{t}^{q_1}L_{x_1}^{q_2}L_{x_2}^{q_3}(\mathbb{R}^{3}) \to
L_{x_1}^{p_1}L_{x_2}^{p_2}L_t^{p_3}(\mathbb{R}^{3})$ is a well
restriction operator.

\item[\rm (4)] If $\wedge^3_{i=1} p_i > (\vee^3_{i=1} q_i) \vee ( q_1q_3 /q_2)$,
     then $T: L_{x_1}^{q_1}L_{x_2}^{q_2}L_t^{q_3}(\mathbb{R}^{3})
     \to
     L_{x_2}^{p_2} L_{x_1}^{p_1}L_t^{p_3} (\mathbb{R}^{3})$
     is a well restriction operator.
     \end{itemize}
\end{prop}

Let $f \in L_{x_1}^{q_1}L_{x_2}^{q_2}L_t^{q_3}(\mathbb{R}^{3})$ so
that
$\|f\|_{L_{x_1}^{q_1}L_{x_2}^{q_2}L_t^{q_3}(\mathbb{R}^{3})}=1$.
Define $F: \mathbb{R}\rightarrow [0, 1]$ by
\begin{align}
F(t):=\left \|\left(\int_{-\infty}^t |f(s, x )|^{q_3}ds
\right)^{1/q_3} \right\|_{L_{x_1}^{q_1}L_{x_2}^{q_2}}^{q_1}
\label{A8}
\end{align}

\begin{lem} \label{lemA.2}
Let $I\subset [0, 1]$ is an interval, then it holds:
\begin{align}
\|\chi_{F^{-1}(I)}f\|_{L_{x_1}^{q_1}L_{x_2}^{q_2}L_t^{q_3}(\mathbb{R}\times\mathbb{R}^{2})}\leq
|I|^{\frac{q_2}{q_1q_3}\wedge\frac{1}{q_1}\wedge \frac{1}{q_2}\wedge
\frac{1}{q_3}}\label{A5}
\end{align}
\end{lem}
\noindent{\bf Proof.} For any $I=(A, B)\subset[0, 1]$, there exist
$t_1, t_2 \in \mathbb{R}$ satisfying
\begin{align}
& A=\Big \|\big(\int_{-\infty}^{t_1} |f(s,
x)|^{q_3}ds\big)^{1/q_3}\Big\|_{L_{x_1}^{q_1}L_{x_2}^{q_2}}, \ \
B=\Big \|\big(\int_{-\infty}^{t_2} |f(s,
x)|^{q_3}ds\big)^{1/q_3}\Big\|_{L_{x_1}^{q_1}L_{x_2}^{q_2}}
 \nonumber
\end{align}
and $F^{-1}(I)= (t_1, t_2)$. For $x=(x_1,x_2)$, we define $J(t, x)$
and $E(t, x_1 )$ by:
\begin{align}
& J(t, x)=\Big(\int_{-\infty}^t|f(s, x)|^{q_3}ds\Big)^{1/q_3},  \ \
 E(t, x_1)=\left( \int J(t,x)^{q_2} dx_2
\right)^{1/q_2}. \label{A14}
\end{align}
It is well known that for $a\ge b>0$,
\begin{align}
& r^{a}-s^{a}\leq C (r^{b}-s^{b})(r^{a-b}+s^{a-b}), \ \ \ 0\leq
s\leq r, \label{b1}
\end{align}
and for $0<a\le b$,
 \begin{align}
   r^{a}-s^{a}\leq
(r^{b}-s^{b})^{a/b} , \ \  0\leq s \leq r. \label{b2}
\end{align}
We divide the proof into the following  four cases.

 {\it Case 1.} $q_3\geq q_2\geq q_1$. From \eqref{b1} we have
\begin{align}
\|\chi_{F^{-1}(I)}f(\cdot, x )\|^{q_3}_{L_t^{q_3}} & \lesssim
(J(t_2, x )^{q_2}-J(t_1, x )^{q_2})J(\infty, x)^{q_3-q_2}\label{21}
\end{align}
 Recalling the assumption
$\|f\|_{L_{x_1}^{q_1}L_{x_2}^{q_2}L_t^{q_3}(\mathbb{R}\times\mathbb{R}^{2})}=1$,
by \eqref{21} , \eqref{A14}, \eqref{b1} and H\"older inequality, we
have
\begin{align}
&\int\Big(\int\|\chi_{F^{-1}(I)} f(\cdot,x)\|^{q_2}_{L_t^{q_3}}dx_2\Big)^{\frac{q_1}{q_2}}dx_1\nonumber\\
& \lesssim \int\Big(\int (J(t_2, x )^{q_2}-J(t_1,
x)^{q_2})^{\frac{q_2}{q_3}}J(\infty,
x)^{(q_3-q_2)\frac{q_2}{q_3}}dx_2\Big)^{\frac{q_1}{q_2}}dx_1\nonumber\\
&\leq\int\Big( \Big\|(J(t_2, x )^{q_2}-J(t_1,
x)^{q_2})^{\frac{q_2}{q_3}}\Big\|_{L_{x_2}^{\frac{q_3}{q_2}}}
\Big\|J(\infty,
x)^{(q_3-q_2)\frac{q_2}{q_3}}\Big\|_{L_{x_2}^{1/(1-q_2/q_3)}}\Big)^{\frac{q_1}{q_2}}dx_1\nonumber\\
&=\int\Big( E(t_2, x_1)^{q_2}-E(t_1,
x_1)^{q_2}\Big)^{\frac{q_1}{q_3}}\Big( E(\infty,
x_1)\Big)^{\frac{(q_3-q_2)q_1}{q_3}}dx_1\label{b3}\\
 & \lesssim \int\Big(
E(t_2, x_1)^{q_1}-E(t_1, x_1)^{q_1}\Big)^{\frac{q_1}{q_3}}\Big(
E(\infty, x_1)\Big)^{\frac{(q_2-q_1)q_1}{q_3}} \Big( E(\infty,
x_1)\Big)^{\frac{(q_3-q_2)q_1}{q_3}}dx_1\nonumber\\
&\leq\Big\|\Big(E(t_2, x_1)^{q_1}-E(t_1,
x_1)^{q_1}\Big)^{\frac{q_1}{q_3}}\Big\|_{L_{x_1}^{q_3/q_1}}
\Big\|E(\infty, x_1)
^{\frac{(q_3-q_1)q_1}{q_3}}\Big\|_{L_{x_1}^{1/(1-q_1/q_3)} }   \label{b6}\\
&\leq (F(t_2)-F(t_1))^{\frac{q_1}{q_3}} F(\infty)^{1-q_1/q_3} \leq
|I|^{\frac{q_1}{q_3}}. \label{A20}
\end{align}

 {\it Case 2.} $q_3\geq q_2, q_2< q_1$. From \eqref{b3}and \eqref{b2}, we have
\begin{align}
&\int\Big(\int\|\chi_{F^{-1}(I)} f(\cdot,x)\|^{q_2}_{L_t^{q_3}}dx_2\Big)^{\frac{q_1}{q_2}}dx_1\nonumber\\
& \lesssim \int\Big( E(t_2, x_1)^{q_2}-E(t_1,
x_1)^{q_2}\Big)^{\frac{q_1}{q_3}}\Big( E(\infty,
x_1)\Big)^{\frac{(q_3-q_2)q_1}{q_3}}dx_1\nonumber\\
&\leq\int\Big( E(t_2, x_1)^{q_1}-E(t_1,
x_1)^{q_1}\Big)^{\frac{q_2}{q_3}}\Big( E(\infty,
x_1)\Big)^{\frac{(q_3-q_2)q_1}{q_3}}dx_1\nonumber\\
&\leq\Big\|\Big(E(t_2, x_1)^{q_1}-E(t_1,
x_1)^{q_1}\Big)^{\frac{q_2}{q_3}}\Big\|_{L_{x_1}^{q_3/q_2}}
\Big\|E(\infty, x_1)
^{\frac{(q_3-q_2)q_1}{q_3}}\Big\|_{L_{x_1}^{1/(1-q_2/q_3)} }   \nonumber\\
&\leq (F(t_2)-F(t_1))^{\frac{q_2}{q_3}} F(\infty)^{1-q_2/q_3} \leq
|I|^{\frac{q_2}{q_3}}. \label{A21}
\end{align}

 {\it Case 3.} $q_3 < q_2\leq q_1$. From \eqref{b2}, we have
\begin{align}
\|\chi_{F^{-1}(I)} f(\cdot,x)\|^{q_3}_{L_t^{q_3}} \leq  (J(t_2, x
)^{q_2}-J(t_1, x )^{q_2})^{q_3/q_2} \label{}
\end{align}
Using \eqref{b2} again, we have
\begin{align}
&\int\Big(\int\|\chi_{F^{-1}(I)} f(\cdot,x)\|^{q_2}_{L_t^{q_3}}dx_2\Big)^{\frac{q_1}{q_2}}dx_1\nonumber\\
&\leq\int\Big(\int(J(t_2, x )^{q_2}-J(t_1, x)^{q_2})dx_2\Big)^{\frac{q_1}{q_2}}dx_1 \label{b5}\\
&=\int\Big( E(t_2, x_1)^{q_2}-E(t_1,
x_1)^{q_2}\Big)^{\frac{q_1}{q_2}}dx_1\nonumber\\
& \le \int\Big( E(t_2, x_1)^{q_1}-E(t_1,
x_1)^{q_1}\Big)dx_1\nonumber\\
&=F(t_2)-F(t_1)=|I|. \label{A22}
\end{align}

 {\it Case 4.} $q_3 < q_2, q_2> q_1$.  From \eqref{b5}, \eqref{b1} and H\"older
inequality we have
\begin{align}
&\int\Big(\int\|\chi_{F^{-1}(I)} f(\cdot,x)\|^{q_2}_{L_t^{q_3}}dx_2\Big)^{\frac{q_1}{q_2}}dx_1 \nonumber\\
&\leq\int\Big(\int(J(t_2, x )^{q_2}-J(t_1, x)^{q_2})dx_2\Big)^{\frac{q_1}{q_2}}dx_1\nonumber\\
&=\int\Big( E(t_2, x_1)^{q_2}-E(t_1,
x_1)^{q_2}\Big)^{\frac{q_1}{q_2}}dx_1\nonumber\\
&\leq\int\Big( E(t_2, x_1)^{q_1}-E(t_1,
x_1)^{q_1}\Big)^{\frac{q_1}{q_2}}E(\infty, x_1)^{\frac{q_1(q_2-q_1)}{q_2}}dx_1\nonumber\\
&\leq\big(F(t_2)-F(t_1)\big)^{\frac{q_1}{q_2}}=|I|^{\frac{q_1}{q_2}}.
\label{A24}
\end{align}
From \eqref{A20}, \eqref{A21}, \eqref{A22} and \eqref{A24} we get
\begin{align}
\|\chi_{F^{-1}(I)}f\|_{L_{x_1}^{q_1}L_{x_2}^{q_2}L_t^{q_3}(\mathbb{R}\times\mathbb{R}^{2})}\leq
C |I|^{\frac{1}{q_1}\wedge \frac{1}{q_2}\wedge
\frac{1}{q_3}\wedge\frac{q_2}{q_1q_3}},
\end{align}
 which
yields \eqref{A5}, as desired. $\hfill\Box$

\begin{figure}
\begin{center}
\includegraphics[height=6cm,width=6cm]{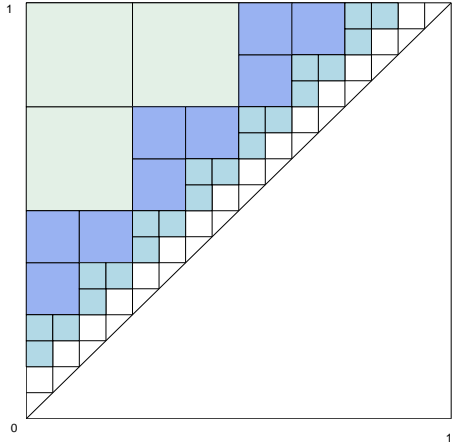}
\begin{minipage}{11cm}
\caption{\small Whitney's decomposition in the triangle.}
\end{minipage}
\end{center}
\end{figure}

\medskip

We will use Whitney's  decomposition to the triangle ${\{(x,
y)\in[0, 1]^2: x<y\}}$ (see Figure 2). First, we divide $[0,1]^2$
into four congruent squares,  consider the square with side-length
$1/2$ in the triangle region and decompose it into four dyadic
squares with side-length $1/4$, then remove the left-upper three
ones in the triangle region. Secondly, considering the remaining
region, we can find three squares with side-length $1/4$ in the
triangle. We decompose each square into four dyadic squares in the
same way as in the first step. Repeating the procedure above to the
end. So, we have decomposed the triangle region into infinite
squares with dyadic border. Let $I$ and $ J$ be the dyadic
subintervals of $[0, 1]$ in the horizontal and perpendicular axes,
respectively. We say that $I\sim J$ if they can consist the
horizontal border and perpendicular border of a square described
above, respectively. From the decomposition above we see that
\begin{enumerate}
\item[(i)]  $|I|=|J|$ and ${\rm dist}(I, \ J)\ge |I| $  for $I\sim J$.
\item[(ii)] The squares in ${\{(x, y)\in[0, 1]^2:
x<y\}}$ are pairwise disjoint.
\item[(iii)] For any dyadic subinterval $J$, there are at most two $I$ with $I
\sim J$.
\end{enumerate}

\noindent{\bf Proof of Proposition \ref{propA.1}.} First, we show
the result of (1). We have
\begin{align}
& T_{re}f(t, x):= \int_{-\infty}^{t}K(t, t')f(t')dt'=\sum_{\{I, J:
I\sim J\}} \chi_{F^{-1}(J)}T(\chi_{F^{-1}(I)}f).
\end{align}
It follows that
\begin{align}
&\|T_{re}f\|_{L_{x_1}^{p_1}L_{x_2}^{p_2}L_t^{p_3}(\mathbb{R}^{3})}
\leq \sum_{j=1}^{\infty} \left\|\sum_{\{I, J: I\sim J, |I|=2^{-j}\}}
\chi_{F^{-1}(J)}T(\chi_{F^{-1}(I)}f)
\right\|_{L_{x_1}^{p_1}L_{x_2}^{p_2}L_t^{p_3}(\mathbb{R}^{3})}.
\label{2}
\end{align}
For any $p\ge 1$, we easily see the following fact:
\begin{align}
& \left\|\sum_{\{I, J: I\sim J, |I|=2^{-j}\}}
\chi_{F^{-1}(J)}T(\chi_{F^{-1}(I)}f)
\right\|^p_{L_t^{p}(\mathbb{R})}
\nonumber\\
& \le 2 \sum_{J_1:  |J_1|=2^{-j} } \int_{\mathbb{R}}
\chi_{F^{-1}(J_1)}|T(\chi_{F^{-1}(J_1)}f) |^p dt. \label{2-a}
\end{align}
Hence, in view of \eqref{2} and \eqref{2-a} we have
\begin{align}
&\|T_{re}f\|_{L_{x_1}^{p_1}L_{x_2}^{p_2}L_t^{p_3}(\mathbb{R}^{3})}
\leq \sum_{j=1}^{\infty} \left\|\left(\sum_{\{I: |I|=2^{-j}\}}
\|T(\chi_{F^{-1}(I)}f)\|^{p_3}_{L_t^{p_3}(\mathbb{R})}
\right)^{1/p_3}\right\|_{L_{x_1}^{p_1}L_{x_2}^{p_2}(\mathbb{R}^{2})}.
\label{2-b}
\end{align}
If $p\le q$, by Minkowski's inequality, we have
\begin{align}
 \left\|\left(\sum_{j}
\|a_j(x,y)\|^{p}_{L_x^p}\right)^{1/p} \right\|_{L^q_y} \le
\left(\sum_{j} \|a_j(x,y)\|^{p}_{L^q_y L_x^p}\right)^{1/p};
\label{2-c}
\end{align}
If $p> q$, in view of $(a+b)^\theta \le a^\theta + b^\theta$ for any
$0\le \theta \le 1, a,b>0$, we have
\begin{align}
 \left\|\left(\sum_{j}
\|a_j(x,y)\|^{p}_{L_x^p}\right)^{1/p} \right\|_{L^q_y} \le
\left(\sum_{j} \|a_j(x,y)\|^{q}_{L^q_y L_x^p}\right)^{1/q}.
\label{2-d}
\end{align}
We divide our discussion into the following three cases.

 {\it Case
1.} $p_1, p_2 \ge p_3$. By  \eqref{2-b},  using \eqref{2-c} twice,
we have
\begin{align}
&\|T_{re}f\|_{L_{x_1}^{p_1}L_{x_2}^{p_2}L_t^{p_3}(\mathbb{R}^{3})}
\leq \sum_{j=1}^{\infty} \left(\sum_{\{I: |I|=2^{-j}\}}
\|T(\chi_{F^{-1}(I)}f)\|^{p_3}_{L_{x_1}^{p_1}L_{x_2}^{p_2}L_t^{p_3}(\mathbb{R}^3)}
\right)^{1/p_3}. \label{2-e}
\end{align}
 {\it Case
2.} $p_1\le p_2 \le p_3$. By  \eqref{2-b}, using \eqref{2-d} twice,
we have
\begin{align}
&\|T_{re}f\|_{L_{x_1}^{p_1}L_{x_2}^{p_2}L_t^{p_3}(\mathbb{R}^{3})}
\leq \sum_{j=1}^{\infty} \left(\sum_{\{I: |I|=2^{-j}\}}
\|T(\chi_{F^{-1}(I)}f)\|^{p_1}_{L_{x_1}^{p_1}L_{x_2}^{p_2}L_t^{p_3}(\mathbb{R}^3)}
\right)^{1/p_1}. \label{2-f}
\end{align}
 {\it Case
3.} $p_2\le p_1 \le p_3$. By \eqref{2-b} and \eqref{2-d}, then
applying \eqref{2-c},  we have
\begin{align}
&\|T_{re}f\|_{L_{x_1}^{p_1}L_{x_2}^{p_2}L_t^{p_3}(\mathbb{R}^{3})}
\leq \sum_{j=1}^{\infty} \left(\sum_{\{I: |I|=2^{-j}\}}
\|T(\chi_{F^{-1}(I)}f)\|^{p_2}_{L_{x_1}^{p_1}L_{x_2}^{p_2}L_t^{p_3}(\mathbb{R}^3)}
\right)^{1/p_2}. \label{2-g}
\end{align}
Denote $p_{\min} = \min(p_1,p_2,p_3)$.  It follows from
\eqref{2-e}--\eqref{2-g} that
\begin{align}
\|T_{re} f\|_{L_{x_1}^{p_1}L_{x_2}^{\infty}L_t^{p_2}} & \lesssim
\sum_{j=1}^{\infty}\left(\sum_{\{I:
|I|=2^{-j}\}}|I|^{\frac{p_{\min}q_2}{q_1q_3}\wedge\frac{p_{\min}}{q_3}\wedge\frac{p_{\min}}{q_2}\wedge
\frac{p_{\min}}{q_1}}\right)^{\frac{1}{p_{\min}}} \nonumber\\
& \lesssim
\sum_{j=1}^{\infty}2^{-j((\frac{q_2}{q_1q_3}\wedge\frac{1}{q_3}\wedge\frac{1}{q_2}\wedge
\frac{1}{q_1})-\frac{1}{p_{\min}})}<\infty.  \label{1.11}
\end{align}
The proof of (4) is almost the same as that of (1) and we omit the
details of the proof.

Next, we prove (2). We have
\begin{align}
\|T_{re}f\|_{L_{t}^{p_1}L_{x_1}^{p_2}L_{x_2}^{p_3}(\mathbb{R}\times\mathbb{R}^{2})}
&\leq \sum_{j=1}^{\infty} \left\|\sum_{\{I, J: I\sim J,
|I|=2^{-j}\}}
\chi_{F^{-1}(J)}T(\chi_{F^{-1}(I)}f) \right\|_{L_{t}^{p_1}L_{x_1}^{p_2}L_{x_2}^{p_3}(\mathbb{R}^{3})}\nonumber\\
&\leq 2\sum_{j=1}^{\infty} \left\|\sum_{\{I:  |I|=2^{-j}\}}
\chi_{F^{-1}(J)}
\|T(\chi_{F^{-1}(I)}f)\|_{L_{x_1}^{p_2}L_{x_2}^{p_3}
(\mathbb{R}^{2})}\right\|_{L^{p_1}_t(\mathbb{R}) }. \nonumber
\end{align}
Using the same way as in \eqref{2-a},
\begin{align}
\|T_{re}f\|_{L_{t}^{p_1}L_{x_1}^{p_2}L_{x_2}^{p_3}(\mathbb{R}^{3})}
& \lesssim \sum_{j=1}^{\infty} \left(\sum_{\{I:  |I|=2^{-j}\}}
\|\chi_{F^{-1}(I)}f\|_{L_{x_1}^{q_1}L_{x_2}^{q_2}L_{t}^{q_3}
(\mathbb{R}^{3})}^{p_1}\right)^{1/p_1}, \nonumber
\end{align}
So, we can control
$\|T_{re}f\|_{L_{t}^{p_1}L_{x_1}^{p_2}L_{x_2}^{p_3}(\mathbb{R}^{3})}$
by the right-hand side of \eqref{1.11} in the case $p_{\min}=p_1$.

Finally, we prove (3). We define $F_1(t)$ as follows.
\begin{align}
F_1(t):= \int_{-\infty}^t \|f(s, x_1,
x_2)\|^{q_1}_{L_{x_1}^{q_2}L_{x_2}^{q_3}}ds. \label{B8}
\end{align}
From the definition of $F_1(t)$, it is easy to see that
\begin{align}
\Big\|\chi_{F^{-1}_1(I)}(s)f(s)\Big\|_{L_{t}^{q_1}L_{x_1}^{q_2}L_{x_2}^{q_3}(\mathbb{R}\times\mathbb{R}^{2})}
= |I|^{1/q_1}. \label{C4}
\end{align}
Hence, replacing \eqref{A5} with \eqref{C4},  we can use the same
way as in the proof of (1) to get the result, as desired. $\hfill
\Box$

\medskip

We can generalized this result to $n$ dimensional spaces:

\begin{lem}
Let $T$ be as in \eqref{A1}. We have the following results.
\begin{itemize}
     \item[\rm (1)] If $\min ( p_1, p_2, p_3) > \max(q_1, q_2, q_3, \ q_1q_3 /q_2)$,
     then $T: L_{x_1}^{q_1}L_{x_2,...,x_n}^{q_2}
L_t^{q_3}(\mathbb{R}^{n+1}) \to L_{x_1}^{p_1}L_{x_2,...,x_n}^{p_2}
L_t^{p_3}(\mathbb{R}^{n+1}) $
     is a well restriction operator.
\item[\rm (2)] If $p_0
> (\vee^3_{i=1} q_i) \vee (q_1q_3/q_2)$, then $T:  L_{x_1}^{q_1}L_{x_2,...,x_n}^{q_2}
L_t^{q_3}(\mathbb{R}^{n+1}) \to L_{t}^{p_0}L_{x_1}^{p_1}...
L_{x_n}^{p_{n}} (\mathbb{R}^{n+1})$ is a well restriction operator.
\item[\rm (3)] If $q_0< \min{(p_1, p_2, p_3)}$, then
$T: L_{t}^{q_0}L_{x_1}^{q_1}... L_{x_n}^{q_n} (\mathbb{R}^{n+1}) \to
L_{x_1}^{p_1}L_{x_2,...,x_n}^{p_2} L_t^{p_3}(\mathbb{R}^{n+1}) $ is
a well restriction operator.
 \item[\rm (4)] If $\min ( p_1, p_2, p_3) > \max(q_1, q_2, q_3, \ q_1q_3 /q_2)$,
     then $T: L_{x_2}^{q_1} L_{x_1,x_3,...,x_n}^{q_2}
L_t^{q_3}(\mathbb{R}^{n+1}) \to L_{x_1}^{p_1}L_{x_2,...,x_n}^{p_2}
L_t^{p_3}(\mathbb{R}^{n+1}) $
     is a well restriction operator.
     \end{itemize}
\end{lem}

\end{appendix}

\noindent{\bf Acknowledgment.}  This work is supported in part by
the National Science Foundation of China, grants  10571004 and
10621061; and the 973 Project Foundation of China, grant
2006CB805902.

\end{document}